\documentclass[a4paper,11pt]{article}
\usepackage[top=30truemm,bottom=30truemm,left=30truemm,right=30truemm]{geometry}
\usepackage{amsmath}
\usepackage{amssymb}
\usepackage{float}
\usepackage{color}
\usepackage{graphicx}
\usepackage{graphics}
\usepackage{graphicx}
\usepackage{tcolorbox}
\usepackage{latexsym}
\usepackage{ascmac}
\usepackage{framed}
\usepackage{fancybox}
\usepackage{theorem}
\usepackage{stmaryrd}
\usepackage{here}
\usepackage{mathtools}
\usepackage{bm}
\usepackage{appendix}
\usepackage{url}

\mathtoolsset{showonlyrefs=true} 

\theoremstyle{break}

\newtheorem{thm}{Theorem}

\makeatletter
\@addtoreset{equation}{section}

\makeatother

\title{{\bf{
Memory-induced blow-up solutions and their dynamical transitions in distributed delay differential equations
}}
}

\author{
Yu Ichida
\thanks{Department of Mathematics, School of Science and Technology, Meiji University, 1-1-1 Higashimita, Tama-ku, Kawasaki 214-8571, Japan, {\tt ichidayu@meiji.ac.jp} }
}

\begin{document}
\maketitle

\begin{abstract}
In this paper, we investigate finite-time blow-up solutions for distributed delay differential equations incorporating memory effects with a specific gamma distribution kernel.
Focusing on typical nonlinear terms that cause finite-time singularities in ordinary differential equations (ODEs), we examine how memory effects change the existence, rate, and qualitative properties of blow-up.
By comparing these behaviors with those of the corresponding non-delayed and memory-free ODEs, we show that time delay originating from memory not only essentially induces finite-time blow-up (``memory-induced blow-up''), but also drastically alters the blow-up profile, such as accelerating algebraic blow-up rates or driving a qualitative transition from ODE quenching to logarithmic blow-up.
Furthermore, by incorporating a self-inhibitory term, we reveal a threshold phenomenon in the phase space that governs the occurrence and non-occurrence of blow-up depending on the initial conditions and parameters.
These results are established by reducing the system to a two-dimensional ODE via the linear chain trick, and analyzing the dynamics at infinity using Poincar\'e-type compactification, blow-up techniques, and the center manifold theorem.
\end{abstract}

{\bf Keywords:
Distributed delay differential equations,
Memory-induced blow-up,
Linear chain trick,
Poincar\'e-type compactification,
Center manifold theorem,
Dynamical transitions
} 

\begin{center}
{\scriptsize Mathematics Subject Classification: 
34C05, 
34K12, 
34K17, 
34K25, 
34K99, 
37C35, 
}
\end{center}

\section{Introduction}
\label{sec:DDDEb-int}
In this paper, we consider the following distributed delay differential equation
\begin{equation}
\dot{u}(t) = f\left( u, \int_{-\infty}^{t} \alpha e^{-\beta(t-s)} g(u(s)) \, ds \right),
\label{eq:DDDEb-int1}
\end{equation}
where $\dot{u}=du/dt$.
The unknown function is $u = u(t)$, where $t$ denotes the present time and $s$ denotes the past time $(-\infty < s \le t)$. 
Let $\alpha$ and $\beta$ be positive constants.
Here, $\alpha$ represents the past history and the strength of memory, $\beta$ is a parameter corresponding to the rate of forgetting history and memory. 
A larger value of $\beta$ implies strong memory decay and a weaker delay effect.
On the other hand, a smaller value of $\beta$ means that the influence of memory and time delay appears strongly. 
Let $\mathbb{R}_{+}$ be the set of positive real numbers, and we impose the initial function
\begin{equation}
u(s)=\phi(s), \quad \phi \in C_{b}:=C_{b}((-\infty, 0]; \mathbb{R}_{+}),
\quad \phi(s)>0,
\quad -\infty < s \le 0,
\label{eq:DDDEb-ini1}
\end{equation}
where $\phi \in C_{b}((-\infty, 0]; \mathbb{R}_{+})$ means that $\phi$ is a positive and continuous function from $(-\infty,0]$ to $\mathbb{R}_{+}$ satisfying the boundedness condition $\sup_{s\le 0} |\phi(s)|<+\infty$.

Furthermore, we define $v=v(t)$ as
\begin{equation}
v(t)=\int_{-\infty}^{t} \alpha e^{-\beta(t-s)} g(u(s)) \, ds.
\label{eq:DDDEb-int2}
\end{equation}
Under the definition of the initial function \eqref{eq:DDDEb-ini1} and $\beta > 0$, the total amount of memory (or history) at the initial time $t = 0$, denoted by $v(0)$, takes a positive finite value $v(0) > 0$.

This paper is to specify $f$ and $g$ as follows.
The goal of this paper is to clarify the behavior of the solutions from the viewpoints of global existence in time and finite-time singularities: 
\begin{quote}
\begin{enumerate}
\item[(F1)] $f(u, v)= v^{p}$ ($p>1$) and $g(u)=u$, that is, \eqref{eq:DDDEb-int1} becomes 
\begin{equation}
\dot{u}(t) = \left(\int_{-\infty}^{t} \alpha e^{-\beta(t-s)}u(s) \, ds \right)^{p}.
\label{eq:DDDEb-int01}
\end{equation}
\item[(F2)] $f(u,v)=(1-v)^{-1}$ and $g(u)=u$, that is, \eqref{eq:DDDEb-int1} leads
\begin{equation}
\dot{u}(t) = \frac{1}{1 - \int_{-\infty}^{t} \alpha e^{-\beta(t-s)} u(s) \, ds}.
\label{eq:DDDEb-int02}
\end{equation}
\item[(F3)] Let $\gamma$ and $a$ be positive constants, with $f(u, v)=(1-v)^{-1}-\gamma u$ and $g(u)=u-a$. 
\eqref{eq:DDDEb-int1} becomes
\begin{equation}
\dot{u}(t) = -\gamma u(t) + \frac{1}{1 - \int_{-\infty}^{t} \alpha(u(s)-a) e^{-\beta(t-s)}  \, ds}.
\label{eq:DDDEb-int03}
\end{equation}
\end{enumerate}
\end{quote}

The function $f$ specified above originates from two typical and classical examples concerning finite-time singularities of solutions to ordinary differential equations (ODEs). 
Here, we introduce two examples of nonlinearities that induce finite-time blow-up and quenching in ODEs.

First, for $p > 1$, consider the following initial value problem:
\begin{equation}
\dot{u} = u^p, \quad u(0) > 0.
\label{eq:DDDEb-int3}
\end{equation}
It is well known that the solution to the initial value problem \eqref{eq:DDDEb-int3} blows up in finite time. 
That is, there exists a finite time $T \in (0, +\infty)$ such that
\begin{equation}
\lim_{t \to T^-} u(t) = +\infty \quad (0 < T < +\infty)
\label{eq:DDDEb-int4}
\end{equation}
with $T^{-}=T-0$.
The time $T$ is called the blow-up time.
Furthermore, the blow-up rate is given by
\begin{equation}
u(t) = C (T-t)^{-\frac{1}{p-1}} (1+o(1)) 
\quad {\rm{as}} \quad t\to T^{-}
\label{eq:DDDEb-int5}
\end{equation}
with a constant $C>0$. 
In this paper, we refer to the algebraic blow-up rate \eqref{eq:DDDEb-int5} as the ODE blow-up rate.

As a second example, consider the initial value problem for the following ODE:
\begin{equation}
\dot{u} = \frac{1}{1 - u}, \quad 0 < u(0) < 1.
\label{eq:DDDEb-int6}
\end{equation}
Then, the solution of \eqref{eq:DDDEb-int6} is given by $u(t) = 1 - \sqrt{(1 - u_0)^2 - 2t}$, and there exists a finite time $T > 0$ such that
\begin{equation}
 \lim_{t \to T^-} u(t) = 1, \quad \lim_{t \to T^-} \vert{}\dot{u}(t)\vert{} = +\infty
\label{eq:DDDEb-int7}
\end{equation}
holds.
In this phenomenon, the time derivative of the solution blows up in finite time, which is referred to as quenching in ODEs.

As can be seen from these simple examples, understanding when and how finite-time singularities occur is crucial. 
For instance, in the context of blow-up solutions, key mathematical issues include the relationship between the existence of blow-up solutions and initial data, the estimation and evaluation of blow-up time, and the derivation of blow-up rates. 
For further details on finite-time singularities such as blow-up and quenching of solutions to differential equations, see, for instance, Matsue \cite{Matsue1, Matsue2} and the references therein.

In real-world phenomena, present state changes do not simply depend on present values; they are influenced by past history and the accumulation of memory, which cannot be ignored. 
Considering time delay to express past history and memory leads to delay differential equations (DDEs). 
Numerous past studies have revealed that time delay possesses a dual functionality, acting either to stabilize or to destabilize solutions (see, e.g., \cite{EIIN2021, EJ2006, GNR2018, IN2021, Ruan2006, Smith, Yagasaki2021} and the references therein).

In this paper, since (F1) and (F2) are imposed, finite-time singularities are expected to occur in the solutions of \eqref{eq:DDDEb-int1}. 
Furthermore, by adding a self-inhibitory term as in (F3), inhibitory behavior is expected. 
Clarifying whether time delay, past history, and forgetting and accumulation of memory induce or suppress finite-time singularities, as well as how thresholds, blow-up mechanisms, and profiles change, leads to the prediction and control of system catastrophes caused by time delay and memory effects, making this work extremely valuable.

Here, we will state several previous studies on finite-time singularities of solutions to delay differential equations.
Ezzinbi-Jazar (\cite{EJ2006}) considered scalar differential equations with a constant delay and clarified the effects of additive and multiplicative delay terms on the induction and suppression of finite-time blow-up by comparing them with ordinary differential equations.
To the best of the author's knowledge, the first result demonstrating that time delay essentially induces blow-up solutions is by Eremin-Ishiwata-Ishiwata-Nakata (\cite{EIIN2021}).
Considering the normal form of the Hopf bifurcation as a non-delayed system, they proved delay-induced blow-up, where solution blow-up occurs by introducing a constant delay.
Yagasaki (\cite{Yagasaki2021}) showed, based on the normal form of the subcritical Hopf bifurcation, that finite-time blow-up of solutions occurs under appropriate conditions using the averaging method.
As a scalar equation with a constant delay investigating the effects of time delay, there is a study by Ishiwata-Nakata (\cite{IN2021}).
However, as far as the author knows, there are almost no results on blow-up solutions arising from distributed delays where past history and forgetting and accumulation of memory.

As can be seen from the way initial functions are given, delay differential equations are infinite-dimensional dynamical systems whose solution behavior must be carefully analyzed. 
Due to the difficulty in handling time delays, analyzing them --including numerical calculations-- is far more challenging than analyzing ordinary differential equations. 
For these reasons, the mechanisms underlying finite-time singularities in solutions of delay differential equations are not yet fully understood.
Therefore, as a first step to directly observe the effects of time delay by focusing on the parameters of memory expressing time delay --while essentially reducing the behavior of solutions to ordinary differential equation systems-- this paper adopts equation \eqref{eq:DDDEb-int1}. 
In addition, by selecting $f$ and $g$ from typical nonlinearities that cause finite-time singularities in ODE solutions, namely \eqref{eq:DDDEb-int3} and \eqref{eq:DDDEb-int6}, the objective of this paper is to investigate how memory effects appearing as time delays essentially affect finite-time singularities. 
We examine this from the viewpoints of the existence or non-existence of finite-time singularities, the rate of singularity formation, acceleration and suppression of blow-up depending on the parameter values and initial memory, and threshold phenomena in blow-up behavior.

In this paper, we set $f$ and $g$ in \eqref{eq:DDDEb-int1} as (F1), (F2), and (F3).
For each case, we investigate the existence or non-existence of finite-time singularities and, when they occur, the profile of singularity formation from the perspective of asymptotic behavior. 
Furthermore, as discussed in Section \ref{sec:DDDEb-d}, we examine how memory accumulation and forgetting relate to blow-up solutions by comparing them with the memory-free non-delayed differential equations obtained as $\beta \to +\infty$. 
We perform discussions for both the case $\alpha = \beta$ and the case where $\alpha \ne \beta$ (i.e., when $\alpha$ and $\beta$ are independent). 
This allows us to determine whether finite-time singularities of solutions are induced by the effects of time delay expressing memory effects.

In this paper, by applying the Linear Chain Trick (see, e.g., \cite{MacDonald1978, MacDonald1989, Ruan2006, Smith} and references therein), Equation \eqref{eq:DDDEb-int1} can be rewritten as a two-dimensional ODEs. 
Unlike differential equations with constant delays, the time-delay information in this system is expressed through the forgetting and accumulation of memory, allowing the problem to be essentially reduced to analyzing a two-dimensional phase space as an invariant set within an infinite-dimensional space.
By applying a combination of Poincar\'e-type compactification (see, e.g., \cite{FAL, Matsue1, Matsue2} and references therein), desingularization of blow-up techniques (see, e.g., \cite{AFJ, MB, FAL} and references therein), and the center manifold theorem (see, e.g., \cite{carr} and references therein) to the resulting two-dimensional system, we fully clarify the qualitative dynamics of this two-dimensional system including at infinity. 
This enables us to obtain detailed information regarding finite-time singularities of the solutions. 
This idea is based on the dynamical systems approach to handling blow-up solutions of ODEs developed in \cite{Matsue1, Matsue2}.

Under assumption (F1), finite-time blow-up solutions exist; while their blow-up rate corresponds to the rate observed in standard ODEs, the speed of blow-up is accelerated. 
On the other hand, under assumption (F2), finite-time blow-up solutions also exist, but the blow-up rate exhibits a logarithmic rate, which differs from the standard ODE rate. 
Furthermore, by imposing (F3) which adds a self-inhibitory term to (F2), we provide results on threshold phenomena and the suppression and control of blow-up (i.e., the existence or non-existence of blow-up solutions depending on parameter values and initial conditions in the two-dimensional system).

This study is the first attempt not only to compare how memory effects expressed as time delays induce or suppress finite-time singularities relative to typical ODE examples, but also to perform a rigorous analysis using a dynamical systems approach. 
In particular, by deriving non-delayed ODEs in which memory effects are absent, we elucidate how time delays resulting from memory effects contribute to solution blow-up, as well as how the qualitative dynamics and transitions of solutions brought about by delays change. 
Consequently, this work serves as the first attempt to explicitly demonstrate ``memory-induced blow-up'' and ``qualitative transitions induced by memory'', offering a fundamentally distinct and novel perspective to the study of finite-time singularities in solutions to delay differential equations.

The rest of this paper is organized as follows.
In the next section, we state the main results of this paper. 
In Sections \ref{sec:DDDEb-pr1}, \ref{sec:DDDEb-pr2}, and \ref{sec:DDDEb-pr3}, we prove Theorems \ref{th:DDDEb-mr1}, \ref{th:DDDEb-mr2}, and \ref{th:DDDEb-mr3}, respectively.
In Section \ref{sec:DDDEb-d}, we derive the non-delayed (memory-free) ordinary differential equation, which removes the memory effect representing time delay, as $\beta \to +\infty$.
Furthermore, by comparing the behavior of the non-delayed equation with the main results of this paper, we demonstrate the occurrence of memory-induced blow-up.

\section{Main results}
\label{sec:DDDEb-mr}
We state the main results of this paper. 
In Theorems \ref{th:DDDEb-mr1} and \ref{th:DDDEb-mr2}, which impose (F1) and (F2), respectively, setting the initial function as in \eqref{eq:DDDEb-ini1} allows the Linear Chain Trick to be effectively applied, providing results on the existence of blow-up solutions and their blow-up rates. 
On the other hand, assuming (F3) reveals a threshold for whether blow-up occurs or not, depending on the initial functions and parameter values.

\begin{thm}
\label{th:DDDEb-mr1}
Let $\phi \in C_b$ be an initial function satisfying $\phi(s) > 0$ for all $s \in (-\infty, 0]$.
For any $\alpha>0$ and $\beta>0$, there exists a finite time $T>0$ such that a solution $u(t)$ of the distributed delay differential equation \eqref{eq:DDDEb-int1} with (F1) blows up at $t=T$.
In addition, the blow-up rates are 
\begin{equation}
u(t) = A_{1} (T-t)^{-\frac{p+1}{p-1}} (1+o(1))
\quad {\rm{as}} \quad t\to T^{-}
\label{eq:DDDEb-mr1}
\end{equation}
and 
\begin{equation}
\dot{u}(t) = A_{2} (T-t)^{-\frac{2p}{p-1}} (1+o(1))
\quad {\rm{as}} \quad t\to T^{-},
\label{eq:DDDEb-mr2}
\end{equation}
where $A_1$ and $A_2$ are positive constants.
\end{thm}

\begin{thm}
\label{th:DDDEb-mr2}
Let $\phi \in C_b$ be an initial function satisfying $\phi(s) > 0$ for all $s \in (-\infty, 0]$.
For any $\alpha>0$ and $\beta>0$, there exists a finite time $T>0$ such that a solution $u(t)$ of the distributed delay differential equation \eqref{eq:DDDEb-int1} with (F2) blows up at $t=T$ with the following rates:
\begin{equation}
u(t) = A_{3} \left(\log\dfrac{1}{T - t}\right)^{1/2} (1+o(1))
\quad {\rm{as}} \quad t\to T^{-}
\label{eq:DDDEb-mr3}
\end{equation}
and 
\begin{equation}
\dot{u}(t) = A_{4}(T-t)^{-1}\left(\log\dfrac{1}{T - t}\right)^{-1/2} (1+o(1))
\quad {\rm{as}} \quad t\to T^{-}.
\label{eq:DDDEb-mr4}
\end{equation}
Here, $A_3$ and $A_4$ are positive constants.
\end{thm}

\begin{thm}
\label{th:DDDEb-mr3}
Assume (F3) for \eqref{eq:DDDEb-int1} holds with positive constants $\gamma$ and $a$.
Let $\phi \in C_b$ be an initial function satisfying $\phi(s) > 0$ for all $s \in (-\infty, 0]$.
Then, there exists a threshold curve in the phase space $(u,w)$ dividing finite-time blow-up and non-blow-up solutions, where $w=(1-v)^{-1}$ and \eqref{eq:DDDEb-int2}.
If the solution blows up in finite time, the blow-up rates are given by \eqref{eq:DDDEb-mr3} and \eqref{eq:DDDEb-mr4}.
\end{thm}

\section{Proof of Theorem \ref{th:DDDEb-mr1}}
\label{sec:DDDEb-pr1}
In the delay differential equation \eqref{eq:DDDEb-int1}, by defining $v=v(t)$ as in \eqref{eq:DDDEb-int2}, we have 
\begin{equation}
\dot{v}= -\beta v(t)+\alpha g(u(t)). 
\label{eq:DDDEb-pr1}
\end{equation}
Thus, \eqref{eq:DDDEb-int1} is reduced to considering the following two-dimensional ODEs.
This method is called the Linear Chain Trick.
For details, see \cite{MacDonald1978, MacDonald1989, Ruan2006, Smith} and the references therein.
Namely, it is sufficient to derive the behavior of solutions, the existence of blow-up solutions, and their blow-up rates for the following two-dimensional ODEs:
\begin{equation}
\begin{cases}
\dot{u}=f(u, v), \\
\dot{v}=  \alpha g(u(t)) -\beta v(t).
\end{cases}
\label{eq:DDDEb-pr2}
\end{equation}
In this section, under assumption (F1), it is sufficient to consider the problem of investigating the behavior of the following system as $u \to +\infty$:
\begin{equation}
\begin{cases}
\dot{u} = v^{p}, \\
\dot{v} = \alpha u-\beta v.
\end{cases}
\label{eq:DDDEb-pr3}
\end{equation}

\subsection{Dynamics on the finite equilibrium for \eqref{eq:DDDEb-pr3}}
\label{sub:DDDEb-pr1}
In \eqref{eq:DDDEb-pr3}, we define the region $\Phi_{1}$ as follows:
\[
\Phi_{1}=\{(u,v) \mid u\ge 0,\,\, v\ge 0\}.
\]
Then, $\Phi_{1}$ is an invariant region.
That is, on the $(u,v)$-plane, trajectories originating in $\Phi_1$ cannot leave it across the axes. 
System \eqref{eq:DDDEb-pr3} has the finite equilibrium $E_{0}:(u,v)=(0,0)$.
The Jacobian matrix $J(E_{0})$ of the vector field \eqref{eq:DDDEb-pr3} in the equilibrium $E_{0}$ is 
\[
J(E_{0})=\left(\begin{array}{cc}
0 & 0 \\ \alpha & -\beta
\end{array}\right).
\]
The eigenvalues of $J(E_{0})$ are $0$ and $-\beta$, and the corresponding eigenvectors are 
\[
{\mathbf{v}}_{1}=(\beta, \alpha)^{T}, \quad 
{\mathbf{v}}_{2}=(0, 1)^{T}
\]
with $T$ as the transpose matrix.
Since the finite equilibrium $E_{0}$ is not hyperbolic, we use the normal form transformation and the center manifold theorem to study the dynamics near $E_{0}$ (for instance, see \cite{carr, FKPM} and references therein).
Let $P$ be a matrix $P=(\mathbf{v}_{1},\mathbf{v}_{2})$.
Then we have
\begin{align*}
\left(\begin{array}{cc}
\dot{u} \\
\dot{v}
\end{array}
\right) &= \left(\begin{array}{cc}
0 & 0 \\ \alpha & -\beta
\end{array}
\right)\left(\begin{array}{cc}
u \\
v
\end{array}
\right)+\left(\begin{array}{cc}
v^{p} \\
0
\end{array}
\right) \\
&= P\left(\begin{array}{cc}
0 & 0 \\
0 & -\beta
\end{array}
\right)P^{-1}\left(\begin{array}{cc}
u \\
v
\end{array}
\right)+\left(\begin{array}{cc}
v^{p} \\
0
\end{array}
\right).
\end{align*}
In addition, we set 
$\left(\begin{array}{cc}
\hat{u} \\
\hat{v}
\end{array}
\right)=P^{-1}\left(\begin{array}{cc}
u\\
v
\end{array}
\right)$.
Hence, we then obtain the following system:
\begin{equation}
\begin{cases} 
d\hat{u}/dt
=  \beta^{-1}(\alpha\hat{u}+\hat{v})^{p},
\\ 
d\hat{v}/dt
= -\beta\hat{v}-\alpha\beta^{-1}(\alpha\hat{u}+\hat{v})^{p}.
\end{cases}
\label{eq:DDDEb-pr4}
\end{equation}
The center manifold theorem is applicable to study the dynamics of \eqref{eq:DDDEb-pr4}.
It implies that there exists a function $h(\hat{u})$ satisfying 
\[
h(0)=\dfrac{dh}{d\hat{u}}(0)=0
\]
such that the center manifold of the origin for \eqref{eq:DDDEb-pr4} is locally represented as $\{(\hat{u},\hat{v}) \,|\, \hat{v}(t)=h(\hat{u}(t))\}$.
 Differentiating it with respect to $t$, we obtain that the approximation of the (graph
of) center manifold is
\begin{equation}
\left\{ (\hat{u}(t), \hat{v}(t)) \mid \hat{v} = -\dfrac{\alpha^{p+1}}{\beta^{2}}\hat{u}^{p} + \mathcal{O}(\hat{u}^{p+1}) \right\}.
\label{eq:DDDEb-pr5}
\end{equation}
In addition, the locally dynamics of two-dimensional system \eqref{eq:DDDEb-pr5} can be understood in terms of a dynamical system on a one-dimensional center manifold:
\begin{equation}
\dfrac{d\hat{u}}{dt} = \dfrac{\alpha^{p}}{\beta}\hat{u}^p + \mathcal{O}(\hat{u}^{p+1}).
\label{eq:DDDEb-pr6}
\end{equation}
By using the transformations in this subsection, we conclude that the approximation of the (graph of) center manifold in \eqref{eq:DDDEb-pr3} is
\begin{equation}
\left\{ (u, v) \mid v= \dfrac{\alpha}{\beta}u-\dfrac{\alpha^{p+1}}{\beta^{p+2}}u^{p}+ \mathcal{O}(u^{p+1}) \right\}
\label{eq:DDDEb-pr7}
\end{equation}
 and the locally dynamics of \eqref{eq:DDDEb-pr3} near $E_{0}$ can be understood in terms of a dynamical system on a one-dimensional center manifold:
\begin{equation}
\dot{u}= \left(\dfrac{\alpha}{\beta}\right)^{p} u^p + \mathcal{O}(u^{p+1}).
\label{eq:DDDEb-pr8}
\end{equation}
Therefore, we see that there exists no trajectory that is attracted to the origin.

\subsection{Short review of the Poincar\'e-type compactification}
\label{sub:DDDEb-pr2}
Originally, the problem under consideration is the behavior of system \eqref{eq:DDDEb-pr3} as $u \to +\infty$. 
Specifically, by combining the results from the previous section with the Poincar\'e-Bendixson theorem (see, e.g., \cite{Wiggins} and references therein), it can be seen that any trajectory originating from an initial value in the domain $\Phi_1$ has no choice but to head toward $u \to +\infty$.

In this section, we state the short review of the Poincar\'e-type compactification that enables the extraction of the behavior of systems \eqref{eq:DDDEb-pr2} and \eqref{eq:DDDEb-pr3} as $u \to +\infty$.
In next sections, we investigate the dynamics at infinity of system \eqref{eq:DDDEb-pr3} induced by this method on $\Phi_2 = \Phi_1 \cup \{ \Vert{}(u,v)\Vert{} = +\infty \}$. 
In particular, $\Phi_2$ is referred to as the Poincar\'e-type disk in this paper. 
For details, see \cite{FAL, Matsue1, Matsue2} and references therein.

\begin{figure}[t]
\centering
\includegraphics[scale=0.3]{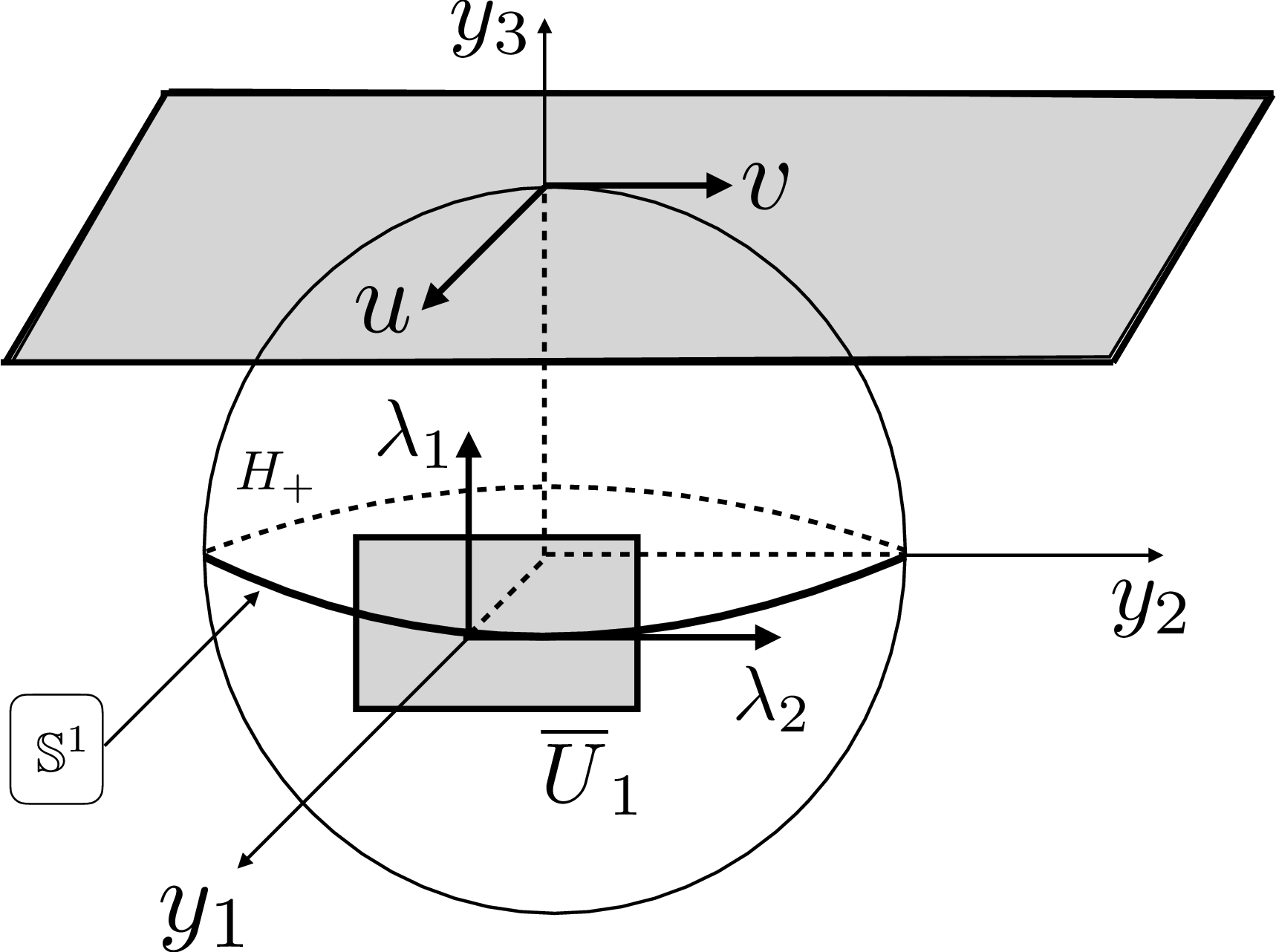}
\caption{Schematic pictures of the locations of the chart $\overline{U}_1$.}
\label{fig:DDDEb1}
\end{figure}

For a coordinate system $(y_1, y_2, y_3) \in \mathbb{R}^3$ and the sphere 
\[
\mathbb{S}^{2}= \{ y \in \mathbb{R}^{3} \, |\, y_{1}^{2} + y_{2}^{2}+y_{3}^{2}=1\} \subset \mathbb{R}^{3},
\]
we embed the phase plane $(u,v)$ of system \eqref{eq:DDDEb-pr2} as $(y_1, y_2, y_3) = (u, v, 1)$. 
Furthermore, we consider the partition of the sphere $\mathbb{S}^2$ given by:
\begin{align*}
\mathbb{S}^{1} &=\{y \in \mathbb{S}^{2}\, | \, y_{3}=0\},\\
H_{+} &= \{ y \in \mathbb{S}^{2}\,|\,y_{3}>0\}, \quad
H_{-} = \{ y \in \mathbb{S}^{2}\,|\,y_{3}<0\}, \\
U_{k} &= \{y \in \mathbb{S}^{2} \, | \, y_{k}>0\}, \quad
V_{k} = \{y \in \mathbb{S}^{2} \, | \, y_{k}<0\}
\end{align*}
with $k=1, 2, 3$.
We define the projections $f^{\pm}$ from $(u, v)\in \mathbb{R}^{2}$ as 
\[
f^{\pm}(u,v):= \pm \left( \dfrac{u}{\Delta(u,v)},\dfrac{v}{\Delta(u,v)},\dfrac{1}{\Delta(u,v)} \right),
\]
where $\Delta(u,v)$ is defined as $\Delta(u,v) = \sqrt{u^{2}+v^{2}+1}$.
Furthermore, we define the local projections $g^{+}_{k} : U_{k} \to \overline{U}_{k}\subset \mathbb{R}^{2}$ and $g^{-}_{k} : V_{k} \to \overline{V}_{k}\subset \mathbb{R}^{2}$ from each $U_k$ and $V_k$ by
\[
g^{+}_{k}(y_{1},y_{2},y_{3}) = - g^{-}_{k}(y_{1},y_{2},y_{3})
 = (y_{m}/y_{k}, y_{n}/y_{k}),
 \]
where $m<n$ with $m,n \neq k$ ($k=1,2,3$), $\overline{U}_{k} = \{y \in \mathbb{R}^{3} \, | \, y_{k} = 1\}$, and $\overline{V}_{k} = \{y \in \mathbb{R}^{3} \, | \, y_{k} = -1\}$.
In addition, denoting the coordinate system on $\overline{U}_{k}$ and $\overline{V}_{k}$ by $(\lambda_{2}, \lambda_{1})$, we have 
\[
(\lambda_{2}, \lambda_{1}) = g^{\pm}_{k}(y).
\]
For instance, when $k=1$, the projection from $\mathbb{R}^2$ to $\overline{U}_1$ is given by 
\[
(g^{+}_{1} \circ f^{+})(u,v) 
= \left ( \dfrac{v}{u},\dfrac{1}{u}\right) = (\lambda_2, \lambda_1). 
\]
See also Figure \ref{fig:DDDEb1}.
Consequently, the dynamical system on $\overline{U}_1$ is obtained via the coordinate transformation $u = 1/\lambda_1$ and $v = \lambda_2/\lambda_1$.

In this paper, since we are interested in the behavior as $u \to +\infty$ within the domain $\Phi_1$, it is sufficient to analyze the dynamics in the local charts $\overline{U}_{1}$ and $\overline{U}_{2}$.

\subsection{Dynamics on the chart $\overline{U}_{1}$}
\label{sub:DDDEb-pr3}
By using the transformation
\[
u=\dfrac{1}{\lambda_1}, \quad 
v=\dfrac{\lambda_2}{\lambda_1},
\]
we examine the dynamics on $\{\lambda_1 = 0\}$ in the $(\lambda_1, \lambda_2)$-plane. 
This enables us to extract a part of the dynamics as $u \to +\infty$. 
Through this transformation, we obtain
\begin{equation}
\begin{cases}
d\lambda_{1}/dt = -\lambda_{1}^{2-p}\lambda_{2}^{2}, \\
d\lambda_{2}/dt = \alpha - \beta\lambda_{2}-\lambda_{1}^{1-p}\lambda_{2}^{p+1}.
\end{cases}
\label{eq:DDDEb-pr9}
\end{equation}
Introducing the time-scale transformation $d\eta/dt=\lambda_{1}^{1-p}$, we have
\begin{equation}
\begin{cases}
d\lambda_{1}/d\eta = -\lambda_{1}\lambda_{2}^{p}, \\
d\lambda_{2}/d\eta = \alpha\lambda_{1}^{p-1}-\beta\lambda_{1}^{p-1}\lambda_{2} -\lambda_{2}^{p+1}.
\label{eq:DDDEb-pr10}
\end{cases}
\end{equation}
The equilibrium of \eqref{eq:DDDEb-pr10} on $\{\lambda_1 = 0\}$ is $E_1: (\lambda_1, \lambda_2) = (0, 0)$. 
The Jacobian matrix $J(E_1)$ at the origin for \eqref{eq:DDDEb-pr10} is
\[
J(E_{1})=\left(\begin{array}{cc}
0 & 0 \\ 0 & 0
\end{array}\right),
\]
which means that the origin is a non-hyperbolic equilibrium since the eigenvalues of $J(E_1)$ are double zero eigenvalues. 
The desingularization of vector fields by the blow-up technique is an effective method to study the behavior near its equilibrium (see \cite{AFJ, MB, FAL} and references therein).
That is, we desingularize the origin by introducing the following blow-up coordinates:
\[
\lambda_{1}=r^{p+1}\bar{\lambda}_{1}, \quad 
\lambda_{2}=r^{p-1}\bar{\lambda}_{2}.
\]
Since we are interested in the dynamics on $\Phi_1$, it is sufficient to examine the region satisfying $\lambda_1 \ge 0$ and $\lambda_2 \ge 0$.

\subsubsection{Dynamics on the chart $\{\bar{\lambda}_1=1\}$}
\label{sub:DDDEb-pr31}
By using the transformation $\lambda_{1}=r^{p+1}$, $\lambda_{2}=r^{p-1}\bar{\lambda}_{2}$ and the time-scale transformation $ds/d\eta=r^{p^{2}-p}$,
\begin{equation}
\begin{cases}
dr/ds = -(p+1)^{-1}r\bar{\lambda}_{2}^{p}, \\
d\bar{\lambda}_{2}/ds = -2(p+1)^{-1}\bar{\lambda}_{2}^{p+1} +\alpha -\beta r^{p-1}\bar{\lambda}_{2}
\end{cases}
\label{eq:DDDEb-pr11}
\end{equation}
holds.
The equilibrium of \eqref{eq:DDDEb-pr11} on $\{r=0\}$ is
\[
E_{2}: (r, \bar{\lambda}_{2})= \left(0, \left[2^{-1}\alpha(p+1)\right]^{\frac{1}{p+1}}\right),
\] 
and the linearized matrix $J(E_2)$ near this equilibrium is
\[
J(E_{2})=\left(\begin{array}{cc}
-(p+1)^{-1}\left[2^{-1}\alpha(p+1)\right]^{\frac{1}{p+1}} & 0 \\
0 & -2\left[2^{-1}\alpha(p+1)\right]^{\frac{1}{p+1}}
\end{array}\right).
\]
Hence, the equilibrium $E_2$ is asymptotically stable.
The approximations of the solution near $E_2$ are as follows:
\begin{equation}
\begin{cases}
r(s) = C_{1}e^{-(p+1)^{-1}\left[2^{-1}\alpha(p+1)\right]^{\frac{1}{p+1}}s}(1+o(1)), \\
\bar{\lambda}_{2}(s) =C_{2} e^{-2\left[2^{-1}\alpha(p+1)\right]^{\frac{1}{p+1}}s}(1+o(1)),
\end{cases}
\quad {\rm{as}} \quad s\to +\infty,
\label{eq:DDDEb-pr12}
\end{equation}
where $C_j$ ($j=1,2$) are positive constants depending on the parameters $(\alpha, \beta, p)$ and initial functions.

\subsubsection{Dynamics on the chart $\{\bar{\lambda}_2=1\}$}
\label{sub:DDDEb-pr32}
The  transformation $\lambda_{1}=r^{p+1}\bar{\lambda}_{1}$, $\lambda_{2}=r^{p-1}$, and the time-scale rescaling $ds/d\eta=r^{p^{2}-p}$ yield
\begin{equation}
\begin{cases}
dr/ds = \alpha(p-1)^{-1}r\bar{\lambda}_{1}^{p-1}-\beta(p-1)^{-1}r^{p}\bar{\lambda}_{1} -(p-1)^{-1}r, \\
d\bar{\lambda}_{1}/ds = -(p+1)(p-1)^{-1}\alpha \bar{\lambda}_{1}^{p}+(p+1)(p-1)^{-1}\beta r^{p-1}\bar{\lambda}_{1}^{2}+2(p-1)^{-1}\bar{\lambda}_{1}.
\end{cases}
\label{eq:DDDEb-pr13}
\end{equation}
The equilibria \eqref{eq:DDDEb-pr13} on $\{r = 0\}$ are
\[
E_{3}: (r, \bar{\lambda}_{1}) =(0,0), \quad 
E_{4}\left(0, \left[2\alpha^{-1}(p+1)^{-1}\right]^{\frac{1}{p-1}}\right).
\]
Calculating the corresponding linearized matrices and their eigenvalues reveals that $E_3$ is a saddle and $E_4$ is asymptotically stable.

\begin{figure}[t]
\begin{center}
\includegraphics[scale=0.37]{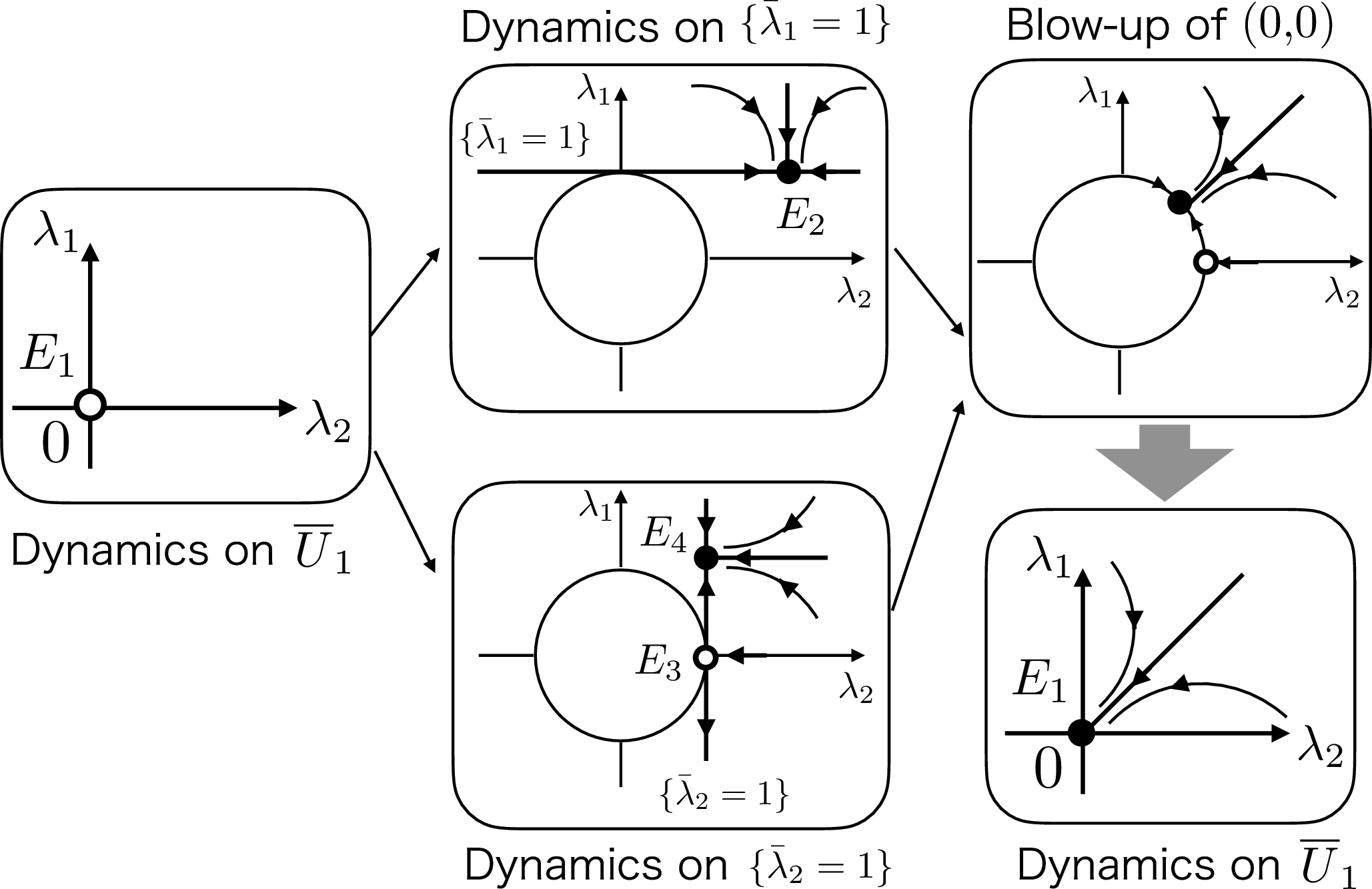}
\caption{Schematic pictures of the dynamics of the blow-up vector fields and $\overline{U}_1$.}
\label{fig:DDDEb2}
\end{center}
\end{figure}

From the above analysis, by combining the behaviors on local charts $\{\bar{\lambda}_1 = 1\}$ and $\{\bar{\lambda}_2 = 1\}$ via desingularization, we can clarify the dynamics near $E_1$. 
See Figure \ref{fig:DDDEb2}.

\subsection{Dynamics on the chart $\overline{U}_{2}$}
\label{sub:DDDEb-pr4}
Using
\[
u=\dfrac{\lambda_2}{\lambda_1}, \quad 
v=\dfrac{1}{\lambda_1}
\]
and $d\eta/dt=\lambda_{1}^{1-p}$, we have 
\begin{equation}
\begin{cases}
d\lambda_{1}/d\eta = -\alpha\lambda_{1}^{p}\lambda_{2}+\beta \lambda_{1}^{p},\\
d\lambda_{2}/d\eta = 1-\alpha\lambda_{1}^{p-1}\lambda_{2}^{2}+\beta \lambda_{1}^{p-1}\lambda_{2}.
\end{cases}
\label{eq:DDDEb-pr14}
\end{equation}
Therefore, there are no equilibrium points, and $d\lambda_2/d\eta > 0$ holds on $\{\lambda_1 = 0\}$.

\subsection{Dynamics and connecting orbits on the Poincar\'e-type disk}
\label{sub:DDDEb-pr5}
\begin{figure}[t]
\begin{center}
\includegraphics[scale=0.25]{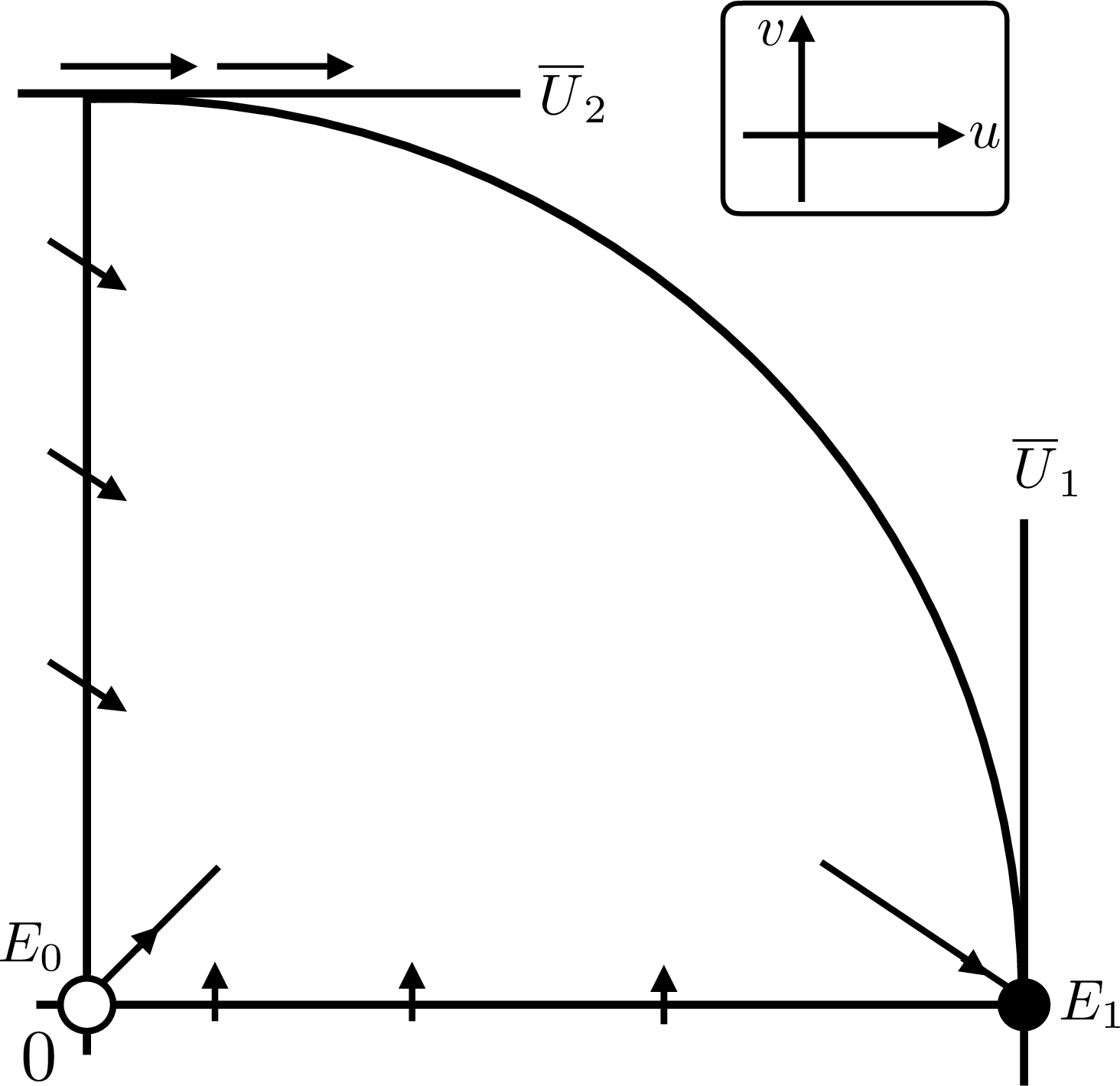}
\caption{Schematic pictures of the dynamics on a Poincar\'e-type disk in \eqref{eq:DDDEb-pr3}.}
\label{fig:DDDEb3}
\end{center}
\end{figure}

By combining the invariant property of the domain $\Phi_1$ defined in Subsection \ref{sub:DDDEb-pr1}, the dynamics near $E_0$, and the dynamics on local charts $\overline{U}_1$ and $\overline{U}_2$, we obtain the dynamics at infinity for \eqref{eq:DDDEb-pr3} (see Figure  \ref{fig:DDDEb3}).

For any compact subset $W \subset \Phi_1$, there are no equilibria or periodic orbits in $W$. 
Thus, by the Poincar\'e-Bendixson theorem, any trajectory starting from an initial point in this subset $W$ cannot stay in $W$ as $t$ increases. 
Namely, any trajectory starting from an arbitrary initial value in $\Phi_1$ has no choice but to be attracted to the equilibrium $E_1$ at infinity. 
This implies that for system \eqref{eq:DDDEb-pr3}, imposing an initial value with $u > 0$ leads to $u \to +\infty$.

\subsection{Proof of Theorem \ref{th:DDDEb-mr1}}
\label{sub:DDDEb-pr6}
Under assumption (F1) in \eqref{eq:DDDEb-int1}, the existence of trajectories in \eqref{eq:DDDEb-int01} such that $u \to +\infty$ has been demonstrated in the preceding subsections. 
In this subsection, we show that the solution blow-up occurs in finite time and derive the blow-up rates \eqref{eq:DDDEb-mr1} and \eqref{eq:DDDEb-mr2}.

Using \eqref{eq:DDDEb-pr12} and the time scale transformations up to the previous subsections, we have
\begin{align*}
\dfrac{ds}{dt}
&=\dfrac{ds}{d\eta}\dfrac{d\eta}{dt}=r^{p^{2}-p}\lambda_{1}^{1-p}=r^{1-p} \\
&\sim \left\{ C_{1}e^{-(p+1)^{-1}\left[2^{-1}\alpha(p+1)\right]^{\frac{1}{p+1}}s}(1+o(1)) \right\}^{1-p}
\quad {\rm{as}} \quad s\to +\infty
\\
&\sim C_{3}e^{\frac{p-1}{p+1}\left[2^{-1}\alpha(p+1)\right]^{\frac{1}{p+1}}s}
\quad {\rm{as}} \quad s\to +\infty,
\end{align*}
where $C_3 > 0$ is a constant, and $F(s) \sim G(s)$ as $s \to +\infty$ means $\lim_{s \to +\infty} |F(s)/G(s)| = 1$. 
Therefore, similarly to \cite{Matsue1, Matsue2}, the maximal existence time $T$ is given by
\[
T:= C_{3} \int_{0}^{+\infty} e^{\frac{p-1}{p+1}\left[2^{-1}\alpha(p+1)\right]^{\frac{1}{p+1}}s}\, ds <+\infty
\]
and the asymptotic relationship is obtained by 
\[
T-t \sim C_{4}e^{-\frac{p-1}{p+1}\left[2^{-1}\alpha(p+1)\right]^{\frac{1}{p+1}}s}
\quad {\rm{as}} \quad s\to +\infty.
\]
Consequently, we obtain the blow-up rate
\begin{align*}
u(t)
&= \lambda_{1}^{-1} =r^{-(p+1)} \\
&\sim \left\{ C_{1}e^{-(p+1)^{-1}\left[2^{-1}\alpha(p+1)\right]^{\frac{1}{p+1}}s}(1+o(1)) \right\}^{-(p+1)}
\quad {\rm{as}} \quad s\to +\infty
\\
&\sim C_{1}e^{\left[2^{-1}\alpha(p+1)\right]^{\frac{1}{p+1}}s}
\quad {\rm{as}} \quad s\to +\infty
\\
&\sim C_{5}(T-t)^{-\frac{p+1}{p-1}}
\quad {\rm{as}} \quad t\to T^{-}.
\end{align*}
By a similar argument, we obtain
\[
v(t) \sim C_{6}(T-t)^{-\frac{2}{p-1}}
\quad {\rm{as}} \quad t\to T^{-}
\]
with positive constants $C_{j}$ ($j=4,5,6$).
Thus, Theorem \ref{th:DDDEb-mr1} is proved.

\section{Proof of Theorem \ref{th:DDDEb-mr2}}
\label{sec:DDDEb-pr2}
We impose (F2) on the delay differential equation \eqref{eq:DDDEb-int1}. 
The delay differential equation under consideration is \eqref{eq:DDDEb-int02}, and by defining $v = v(t)$ as in \eqref{eq:DDDEb-int2}, it is sufficient to consider the behavior of the following two-dimensional system:
\begin{equation}
\begin{cases}
\dot{u} = (1-v)^{-1}, \\
\dot{v} = \alpha u-\beta v.
\end{cases}
\label{eq:DDDEb-pr15}
\end{equation}
Focusing on the range $0 \le v(t) < 1$, we introduce the transformation:
\begin{equation}
w(t)=\dfrac{1}{1-v(t)}.
\label{eq:DDDEb-pr16}
\end{equation}
Hence, the problem is reduced to considering the following system:
\begin{equation}
\begin{cases}
\dot{u} = w, \\
\dot{w} = \alpha uw^{2}-\beta w^{2}+\beta w.
\end{cases}
\label{eq:DDDEb-pr17}
\end{equation}
In \eqref{eq:DDDEb-pr17}, the coordinate axes are invariant sets. 
We define the region $\Phi_3$ as follows:
\[
\Phi_{3}=\{(u,w) \mid u\ge0,\,\, w\ge 1\}.
\]
In particular, since we are interested in the behavior $w \to +\infty$ corresponding to $v \to 1 - 0$, we extract and analyze the behavior as $w \to +\infty$ using Poincar\'e-type compactification in the same manner as in Section \ref{sec:DDDEb-pr1}.

\subsection{Dynamics on the chart $\overline{U}_2$}
\label{sub:DDDEb-pr7}
To investigate the behavior as $w \to +\infty$, we use the transformation
\[
u=\dfrac{\lambda_2}{\lambda_1}, \quad 
w=\dfrac{1}{\lambda_1}
\]
and the time scale transformation $d\eta/dt = \lambda_1^{-2}$, which yields
\begin{equation}
\begin{cases}
d\lambda_{1}/d\eta
= -\alpha\lambda_{1}\lambda_{2}+\beta\lambda_{1}^{2}-\beta\lambda_{1}^{3},
\\
d\lambda_{2}/d\eta
=\lambda_{1}^{2}-\alpha\lambda_{2}^{2}+\beta\lambda_{1}\lambda_{2} -\beta\lambda_{1}^{2}\lambda_{2}.
\end{cases}
\label{eq:DDDEb-pr18}
\end{equation}
The equilibrium of \eqref{eq:DDDEb-pr18} on $\{\lambda_1 = 0\}$ is $e_1: (\lambda_1, \lambda_2) = (0, 0)$. 
Since the corresponding linearized matrix has double zero eigenvalues, we desingularize the point using the transformation:
\[
\lambda_{1}=r\bar{\lambda}_{1}, \quad 
\lambda_{2}=r\bar{\lambda}_{2}.
\]

\subsubsection{Dynamics on the chart $\{\bar{\lambda}_{1}=1\}$}
\label{sub:DDDEb-pre71}
By using the transformation $\lambda_1 = r$, $\lambda_2 = r \bar{\lambda}_2$ and the time scale transformation $ds/d\eta = r$, we obtain
\begin{equation}
\begin{cases}
dr/ds= -\alpha r\bar{\lambda}_{2}+\beta r - \beta r^{2}, \\
d\bar{\lambda}_{2}/ds= 1.
\end{cases}
\label{eq:DDDEb-pr19}
\end{equation}
There are no equilibria on $\{r = 0\}$. 
On $\{r = 0\}$, $d\bar{\lambda}_2/ds > 0$ holds.

\subsubsection{Dynamics on the chart $\{\bar{\lambda}_{2}=1\}$}
\label{sub:DDDEb-pre72}
Transformation $\lambda_{1}=r\bar{\lambda}_{1}$, $\lambda_{2}=r$, and the time-scale transformation are yielded the following system:
\begin{equation}
\begin{cases}
dr/ds= r\bar{\lambda}_{1}^{2} - \alpha r + \beta r\bar{\lambda}_{1} -\beta r^{2}\bar{\lambda}_{1}^{2}, \\
d\bar{\lambda}_{1}/ds = - \bar{\lambda}_{1}^{3}.
\end{cases}
\label{eq:DDDEb-pr20}
\end{equation}
The equilibrium on $\{r=0\}$ is $e_{2}: (r, \bar{\lambda}_{1})=(0,0)$.
The corresponding linearized matrix $J(e_2)$ is
\[
J(e_{2})=\left(\begin{array}{cc}
-\alpha & 0 \\ 0 & 0 
\end{array}\right),
\]
and the dynamics near the origin can be understood via the center manifold theorem. 
The invariant set $\{r = 0\}$ is obtained as one representation of the center manifold, and the dynamics on this center manifold is given by
\begin{equation}
\dfrac{d\bar{\lambda}_{1}}{ds}=-\bar{\lambda}_{1}^{3}.
\label{eq:DDDEb-pr21}
\end{equation}
Note that in general, when $\{r = 0\}$ is an invariant set near $e_2$, the uniqueness of the center manifold is not guaranteed. 
However, since only the leading term of the center manifold is required in this paper and the behavior on any center manifold possesses this leading term, it is sufficient to consider the dynamics on the invariant manifold $r = 0$.

\subsection{Dynamics on the chart $\overline{U}_1$}
\label{sub:DDDEb-pr8}
In Subsection \ref{sub:DDDEb-pr7}, we extracted the behavior as $w \to +\infty$. 
In this subsection, by examining the dynamics on another local chart $\overline{U}_1$, we investigate the behavior at infinity other than that in the previous subsection. 
By using the transformation
\[
u=\dfrac{1}{\lambda_1}, \quad 
w=\dfrac{\lambda_2}{\lambda_1}
\]
and the time scale transformation $d\eta/dt = \lambda_1^{-2}$, we obtain
\begin{equation}
\begin{cases}
d\lambda_{1}/d\eta= -\lambda_{1}^{3}\lambda_{2},
\\
d\lambda_{2}/d\eta= \alpha \lambda_{2}^{2}-\beta\lambda_{1}\lambda_{2}^{2}+\beta \lambda_{1}^{2}\lambda_{2}-\lambda_{1}^{2}\lambda_{2}^{2}.
\end{cases}
\label{eq:DDDEb-pr22}
\end{equation}
The equilibrium of this system on $\{\lambda_{1}=0\}$ is $e_{3}: (\lambda_{1}, \lambda_{2})=(0,0)$, and the corresponding linearized matrix $J(e_{3})$ has double zero eigenvalues.
Thus, as in the previous subsections, we can clarify the dynamics near $e_3$ via the blow-up transformation:
\[
\lambda_{1}=r\bar{\lambda}_{1}, \quad \lambda_{2}=r^{2}\bar{\lambda}_{2}.
\]

\subsubsection{Dynamics on the chart $\{\bar{\lambda}_{1}=1\}$}
\label{sub:DDDEb-pre81}
Using $\lambda_{1}=r$, $\lambda_{2}=r^{2}\bar{\lambda}_{2}$, and the time-rescaling $ds/d\eta=r^{2}$, we obtain 
\begin{equation}
\begin{cases}
dr/ds= -r^{3}\bar{\lambda}_{2}, \\
d\bar{\lambda}_{2}/ds= r^{2}\bar{\lambda}_{2}^{2}+\alpha \bar{\lambda}_{2}^{2}-\beta r\bar{\lambda}_{2}^{2}.
\end{cases}
\label{eq:DDDEb-pr23}
\end{equation}
As an equilibrium on $\{r = 0\}$, $e_4: (r, \bar{\lambda}_2) = (0, 0)$ exists. 
Although its linearized matrix has double zero eigenvalues, $dr/ds = 0$ and $d\bar{\lambda}_2/ds > 0$ hold near $r = 0$, indicating that no trajectory converges to the origin.

\subsubsection{Dynamics on the chart $\{\bar{\lambda}_{2}=1\}$}
\label{sub:DDDEb-pre82}
We introduce the formulas $\lambda_{1}=r\bar{\lambda}_{1}$, $\lambda_{2}=r^{2}$ and time-rescaling $ds/d\eta=r^{2}$.
We then have 
\begin{equation}
\begin{cases}
dr/ds= 2^{-1}\alpha r - 2^{-1}\beta r^{2}\bar{\lambda}_{1} + 2^{-1}\beta r\bar{\lambda}_{1}^{2}-2^{-1}r^{3}\bar{\lambda}_{1}^{2}, \\
d\bar{\lambda}_{1}/ds = -2^{-1}\alpha \bar{\lambda}_{1}+2^{-1}\beta \bar{\lambda}_{1}^{2} -2^{-1}\beta \bar{\lambda}_{1}^{3}-2^{-1}r^{2}\bar{\lambda}_{1}^{3}.
\end{cases}
\label{eq:DDDEb-pr24}
\end{equation}
The equilibrium on $\{r=0\}$ is $e_{5}:(r, \bar{\lambda}_{1})=(0,0)$.
The corresponding linearized matrix $J(e_{5})$ is 
\[
J(e_{5})=\left(\begin{array}{cc}
2^{-1}\alpha & 0 \\ 0 & -2^{-1}\alpha
\end{array}\right),
\]
showing that $e_5$ is a saddle. 
From the above, we see that there is no trajectory attracted to the origin on the local chart $\overline{U}_1$.

\subsection{Proof of Theorem \ref{th:DDDEb-mr2}}
\label{sub:DDDEb-pr9}
\begin{figure}[h]
\begin{center}
\includegraphics[scale=0.21]{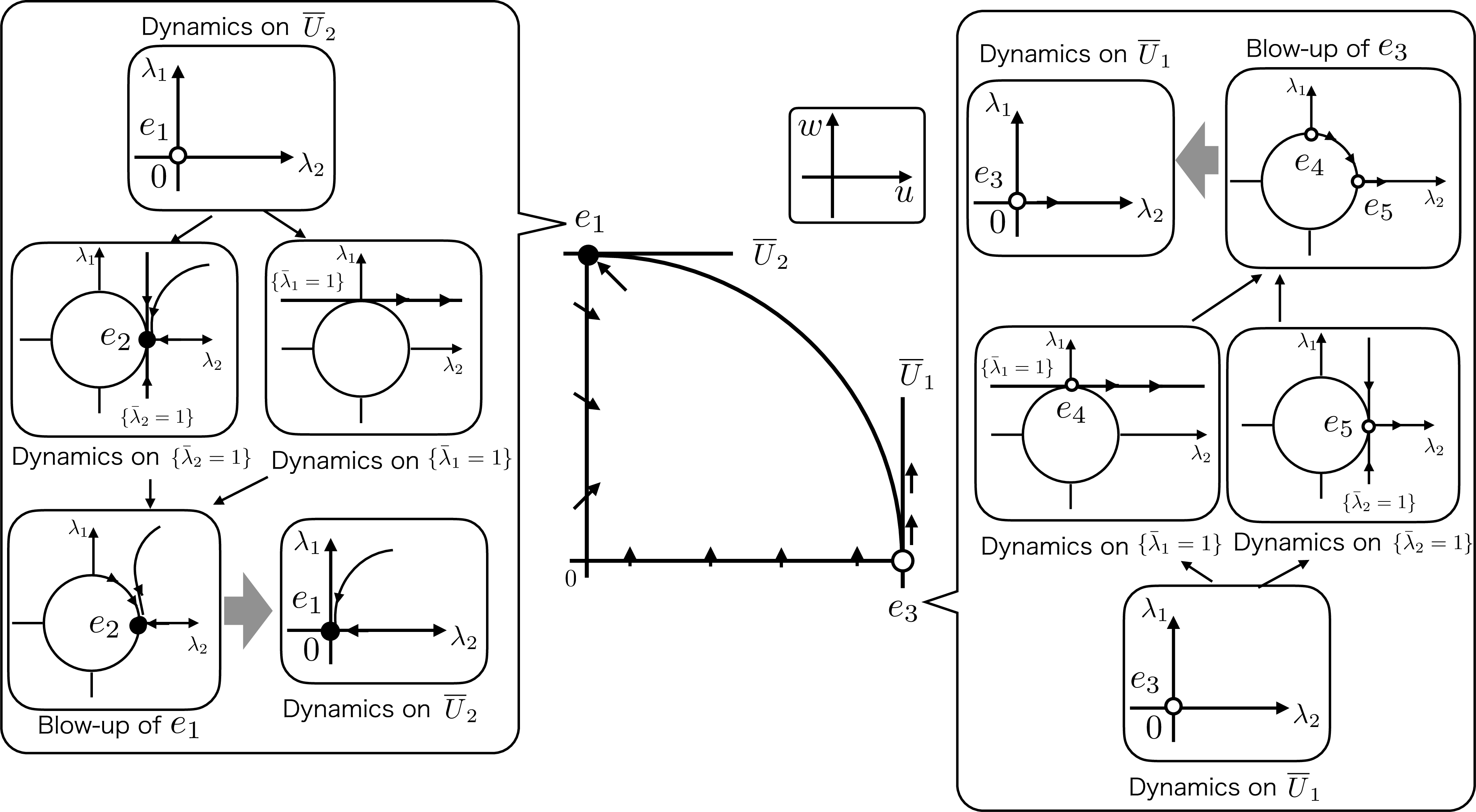}
\caption{Schematic pictures of the dynamics on a Poincar\'e-type disk in \eqref{eq:DDDEb-pr17}.}
\label{fig:DDDEb4}
\end{center}
\end{figure}

By combining the dynamics on local charts $\overline{U}_1$ and $\overline{U}_2$ for \eqref{eq:DDDEb-pr17}, we obtain the dynamics including infinity as shown in Figure \ref{fig:DDDEb4}. 
Since there are no equilibria or periodic orbits on any compact subset in the domain, a similar argument as in Subsection \ref{sub:DDDEb-pr5} shows that any trajectory starting from an initial value in $\Phi_3$ is attracted to the equilibrium $e_1$. 
This means that in \eqref{eq:DDDEb-pr17}, for any initial value in $\Phi_3$ and any parameters $\alpha > 0$ and $\beta > 0$, $w \to +\infty$ holds, which implies $v \to 1 - 0$.

What remains to be shown is that the divergence of $w$ occurs in finite time, along with the derivation of the blow-up rates \eqref{eq:DDDEb-mr3} and \eqref{eq:DDDEb-mr4}. 
The principal part of the dynamics near the equilibrium $e_2$ in Subsection \ref{sub:DDDEb-pre72} is given by $dr/ds \sim -\alpha r$, so with a constant $C_1$, we have
\begin{equation}
r(s) \sim C_{1}e^{-\alpha s}
\quad {\rm{as}} \quad s\to +\infty
\label{eq:DDDEb-as1}
\end{equation}
and the dynamics on the center manifold yields
\begin{equation}
\bar{\lambda}_{1}(s) \sim C_{2}s^{-\frac{1}{2}}
\quad {\rm{as}} \quad s\to +\infty,
\label{eq:DDDEb-as2}
\end{equation}
where $C_2$ is a constant. 
From the relation of time scale transformations, we obtain
\begin{equation}
\dfrac{ds}{dt} 
= \dfrac{ds}{d\eta}\dfrac{d\eta}{dt}
= r\cdot \lambda_{1}^{-2} 
\sim C_{3}\cdot s \cdot e^{\alpha s}
\quad {\rm{as}} \quad s\to +\infty,
\label{eq:DDDEb-as3}
\end{equation}
where $C_3$ is a constant. 
Thus, with a constant $C_4$, we get
\begin{equation}
\dfrac{dt}{ds} \sim C_{4}\dfrac{e^{-\alpha s}}{s}
\quad {\rm{as}} \quad s\to +\infty.
\label{eq:DDDEb-as4}
\end{equation}
Integrating both sides of \eqref{eq:DDDEb-as4} over $[s_0, +\infty)$ with $s_0 > 0$, since $1/s \le 1/s_0$ for $0 < s_0 \le s$, taking $s \to +\infty$ yields
\begin{equation}
t(+\infty) -t(s_{0}) 
=\int_{s_{0}}^{+\infty}\dfrac{dt}{ds}\,ds
\sim C_{4} \int_{s_{0}}^{+\infty} \dfrac{e^{-\alpha s}}{s}\, ds 
\le C_{4} \int_{s_{0}}^{+\infty} \dfrac{e^{-\alpha s}}{s_{0}}\, ds
=\dfrac{C_{4}}{\alpha s_{0}}e^{-\alpha s_{0}}.
\label{eq:DDDEb-as5}
\end{equation}
Hence, we see that $t(+\infty) < +\infty$. 
Here, we set $T$ as
\begin{equation}
T:= t(+\infty)=\lim_{s\to +\infty} t(s).
\label{eq:DDDEb-as6}
\end{equation}
For sufficiently large $s$ (i.e., as $s \to +\infty$), the following estimate holds for any positive constant $\varepsilon > 0$:
\begin{equation}
 (C_4 - \varepsilon) \dfrac{e^{-\alpha \sigma}}{\sigma} 
 \le \dfrac{dt}{d\sigma} 
 \le (C_4 + \varepsilon) \dfrac{e^{-\alpha \sigma}}{\sigma}.
\label{eq:DDDEb-as7}
\end{equation}
Since $|T|<+\infty$, integrating \eqref{eq:DDDEb-as7} over $[s, +\infty)$ gives
\begin{equation}
(C_4 - \varepsilon) \int_{s}^{+\infty} \dfrac{e^{-\alpha \sigma}}{\sigma} d\sigma 
\le T - t(s) 
\le (C_4 + \varepsilon) \int_{s}^{+\infty} \dfrac{e^{-\alpha \sigma}}{\sigma} d\sigma. 
\label{eq:DDDEb-as8}
\end{equation}
Here, integration by parts yields
\begin{equation}
 \int_{s}^{+\infty} \dfrac{e^{-\alpha \sigma}}{\sigma} d\sigma 
 = \dfrac{e^{-\alpha s}}{\alpha s} - \int_{s}^{+\infty} \dfrac{e^{-\alpha \sigma}}{\alpha \sigma^2} d\sigma,
\label{eq:DDDEb-as9}
\end{equation}
which leads to the asymptotic relationship:
\begin{equation}
\int_{s}^{+\infty} \dfrac{e^{-\alpha \sigma}}{\sigma} d\sigma 
\sim \dfrac{e^{-\alpha s}}{\alpha s} 
\quad {\rm{as}} \quad s\to +\infty.
\label{eq:DDDEs-as10}
\end{equation}
Therefore, dividing each side of \eqref{eq:DDDEb-as8} by $(\alpha s)^{-1} e^{-\alpha s}$ and taking $s \to +\infty$, we obtain
\begin{equation}
\dfrac{C_4 - \varepsilon}{\alpha} 
\le \liminf_{s \to +\infty} \dfrac{T - t(s)}{\frac{e^{-\alpha s}}{s}}
\le \limsup_{s \to +\infty} \dfrac{T - t(s)}{\frac{e^{-\alpha s}}{s}} 
\le \dfrac{C_4 + \varepsilon}{\alpha}. 
\label{eq:DDDEb-as12}
\end{equation}
Since $\varepsilon>0$ can be chosen arbitrarily small, taking the limit $\varepsilon\to 0$ yields the asymptotic relation
\begin{equation}
 T -t \sim \dfrac{C_{4}}{\alpha}\dfrac{e^{-\alpha s}}{s}
\quad {\rm{as}} \quad s\to +\infty.
\label{eq:DDDEb-as13}
\end{equation}
From \eqref{eq:DDDEb-as13}, 
\begin{equation}
e^{\alpha s} \sim \dfrac{C_{4}}{\alpha \cdot s\cdot (T-t)}
\quad {\rm{as}} \quad s\to +\infty
\label{eq:DDDEb-as14}
\end{equation}
holds.
Furthermore, taking the logarithm of both sides of \eqref{eq:DDDEb-as13}, with a constant $C_5$, we get
\begin{equation}
\log (T-t) \sim \log \dfrac{C_4}{\alpha} + \log e^{-\alpha s} +\log \dfrac{1}{s}
= C_{5} - \alpha s -\log s 
\sim -\alpha s
\quad {\rm{as}} \quad s\to +\infty.
\label{eq:DDDEb-as15}
\end{equation}
Since $s\to +\infty$ corresponds to $t\to T^{-}$, we have 
\begin{equation}
s \sim \dfrac{1}{\alpha}\log(T-t)^{-1} = \dfrac{1}{\alpha}\log \dfrac{1}{T-t}
\quad {\rm{as}} \quad t\to T^{-}.
\label{eq:DDDEb-as16}
\end{equation}
Thus, using \eqref{eq:DDDEb-as1}, \eqref{eq:DDDEb-as2}, \eqref{eq:DDDEb-as14}, and \eqref{eq:DDDEb-as16}, we obtain 
\begin{align*}
\dot{u}(t)
&=w(t) \\ 
&= \dfrac{1}{\lambda_{1}}
= \dfrac{1}{r\bar{\lambda}_{1}}
= r^{-1}\bar{\lambda}_{1}^{-1}
\\
&\sim \left\{ C_{1}e^{-\alpha s} \right\}^{-1}\cdot \left\{ C_{2}s^{-\frac{1}{2}} \right\}^{-1}
= C_{6} \cdot \sqrt{s}\cdot e^{\alpha s}
\quad {\rm{as}} \quad s\to +\infty
\\
&\sim C_{6}\cdot \sqrt{s} \cdot \dfrac{C_{4}}{\alpha\cdot s\cdot (T-t)}
= \dfrac{C_{4}\cdot C_{6}}{\alpha \cdot (T-t)\cdot \sqrt{s}}
\\
&\sim 
C_{7}\cdot \dfrac{1}{T-t} \cdot \left[ \log(T-t)^{-1} \right]^{-\frac{1}{2}}
 \quad {\rm{as}} \quad t\to T^{-}.
\end{align*}
Similarly, we obtain 
\begin{align*}
u(t)
&=\dfrac{\lambda_{2}}{\lambda_{1}}
= \dfrac{r}{r\bar{\lambda}_1}
= \bar{\lambda}_1^{-1}
\sim  \left\{ C_{2}s^{-\frac{1}{2}} \right\}^{-1}
\quad {\rm{as}} \quad s\to +\infty
\\
&= C_{8} \sqrt{s}
\sim C_{9} \sqrt{\log(T-t)^{-1}}
 \quad {\rm{as}} \quad t\to T^{-}.
\end{align*}
From the above, the asymptotic behaviors \eqref{eq:DDDEb-mr3} and \eqref{eq:DDDEb-mr4} are proved. 
These satisfy $\dot{u} \sim w$.
In addition, we have the following relation:
\[
\dfrac{u}{w} 
\sim \dfrac{C_{9}}{C_{7}}(T-t)\cdot \log(T-t)^{-1}
 \quad {\rm{as}} \quad t\to T^{-}.
\]
Setting $X = T - t$ maps $t \to T^{-}$ to $X \to 0$, 
\[
\lim_{t\to T-0}\dfrac{u}{w} 
=\lim_{X\to 0} \dfrac{u}{w} 
=\lim_{X\to 0} \dfrac{C_{9}}{C_{7}} X\log X^{-1}
=-\lim_{X\to 0} \dfrac{C_{9}}{C_{7}} X\log X
=0
\]
holds.
It implies 
\[
\lambda_{2}=\dfrac{u}{w} \to 0 
\quad {\rm{as}} \quad t\to T^{-}
\]
and the trajectory approaches the boundary at infinity, $\lambda_2 = 0$, tangentially.

\section{Proof of Theorem \ref{th:DDDEb-mr3}}
\label{sec:DDDEb-pr3}
In this section, we consider equation \eqref{eq:DDDEb-int03}, which imposes (F3) on \eqref{eq:DDDEb-int1}.
By the Linear Chain Trick based on \eqref{eq:DDDEb-int2}, we analyze the phase plane induced by the following system:
\begin{equation}
\begin{cases}
\dot{u} = \dfrac{1}{1 - v} - \gamma u,\\
\dot{v} = \alpha(u - a) - \beta v.
\end{cases}
\label{eq:DDDEb-pr25}
\end{equation}
Here, introducing the transformation \eqref{eq:DDDEb-pr16} as in Section \ref{sec:DDDEb-pr2} reduces the problem to considering the system:
\begin{equation}
\begin{cases}
\dot{u} = w - \gamma u, \\
\dot{w} = \alpha(u-a)w^2 - \beta w^2 + \beta w.
\end{cases}
\label{eq:DDDEb-pr26}
\end{equation}
Transformation \eqref{eq:DDDEb-pr16} maps $v \to 1 - 0$ to $w \to +\infty$, and in \eqref{eq:DDDEb-pr26}, it suffices to consider the behavior in the region $\Phi_3$ specified in Section \ref{sec:DDDEb-pr2}.

\subsection{The nondimensionalization of \eqref{eq:DDDEb-pr26}}
\label{sub:DDDEb-pr10}
By introducing the nondimensionalization
\begin{equation}
\tau = \gamma t, \quad 
U = \gamma u, \quad 
W = w
\label{eq:DDDEb-pr27}
\end{equation}
the differential equations satisfied by $U = U(\tau)$ and $W = W(\tau)$ become
\begin{equation}
\begin{cases}
U_{\tau} = -U+W, \\
W_{\tau} = (k_{1}U-k_{2})W^{2}+k_{3}W,
\end{cases}
\label{eq:DDDEb-pr28}
\end{equation}
where the constants $k_1$, $k_2$, $k_3$ are positive constants defined by
\begin{equation}
k_1 = \dfrac{\alpha}{\gamma^2}, \quad
k_2 = \dfrac{\alpha a + \beta}{\gamma}, \quad
k_3 = \dfrac{\beta}{\gamma}.
\label{eq:DDDEb-pr29}
\end{equation}
Therefore, we consider the problem of examining the behavior of equation \eqref{eq:DDDEb-pr28} in the following region $\Phi_4$:
\begin{equation}
\Phi_{4}=\{(U, W) \mid U(\tau)\ge 0,\,\, W(\tau) \ge 1\}.
\label{eq:DDDEb-pr30}
\end{equation}
For equation \eqref{eq:DDDEb-pr28}, as in Sections \ref{sec:DDDEb-pr1} and \ref{sec:DDDEb-pr2}, we examine the finite equilibria and the dynamics in their neighborhood, as well as the dynamics at infinity induced by Poincar\'e-type compactification.
In particular, note that since
\[
k_{2}=\dfrac{\alpha a+\beta}{\gamma}>\dfrac{\beta}{\gamma}=k_{3}>0
\]
holds, $W=k_{3}/k_{2}<1$ holds for $W_{\tau}=0$ at $U=0$.

\subsection{Dynamics on the finite equilibria}
\label{sub:DDDEb-pr11}
The equation satisfied by finite equilibria of \eqref{eq:DDDEb-pr28} is
\begin{equation}
W[k_{1}W^{2}-k_{2}W+k_{3}]=0.
\label{eq:DDDEb-pr31}
\end{equation}
The finite equilibria residing in the domain $\Phi_4$ defined in \eqref{eq:DDDEb-pr30} are as follows:
\begin{itemize}
\item When $k_{2}^{2}-4k_{1}k_{3}<0$, none exists.
\item When $k_{2}^{2}-4k_{1}k_{3}=0$, and $k_{2}\ge 2k_{1}$ holds, there exists a unique equilibrium $P_{0}: (U,W)=(k_{2}/2k_{1}, k_{2}/2k_{1})$ in $\Phi_{4}$.
\item When $k_{2}^{2}-4k_{1}k_{3}>0$, $P_{\pm}: (U, W)=(M_{\pm}, M_{\pm})$ exist.
Here, we define $M_{\pm}$ as 
\begin{equation}
M_{\pm} := \dfrac{k_{2}\pm \sqrt{k_{2}^{2}-4k_{1}k_{3}}}{2k_{1}}>0.
\label{eq:DDDEb-pr32}
\end{equation}
We derive the condition for these equilibria to exist in $\Phi_4$. 
By finding the discriminant, axis equation, and endpoint conditions for the quadratic function $F(W)=k_{1}W^{2}-k_{2}W+k_{3}$ originating from \eqref{eq:DDDEb-pr31}, the necessary and sufficient condition for $M_{\pm}\ge1$ is determined to be $2 k_1 < k_2 \le k_1 + k_3$ and $k_2^2 > 4 k_1 k_3$.
That is,
\begin{enumerate}
\item[(i)] When $2\alpha \le \gamma(\alpha a+\beta) \le \alpha+\gamma\beta$, the finite equilibria $P_{\pm}: (U, W)=(M_{\pm}, M_{\pm})$ exist in $\Phi_{4}$.
\item[(ii)] Otherwise, no finite equilibrium point exists in $\Phi_4$. 
\end{enumerate}
\end{itemize}
The linearized matrices near these finite equilibria are
\begin{align*}
J(P_{0})=\left(\begin{array}{cc}
-1 & 1 \\ k_{3} & -k_{3}
\end{array}\right), \quad
J(P_{\pm})=\left(\begin{array}{cc}
-1 & 1 \\ k_{1}\{M_{\pm}\}^{2} & -k_{3}
\end{array}\right).
\end{align*}
The eigenvalues $\lambda_\pm$ of the linearized matrix $J(P_\pm)$ are
\begin{equation}
\lambda_{\pm} = \dfrac{-(1 + k_3) \pm \sqrt{(1 - k_3)^2 + 4k_1 (M_\pm)^2}}{2},
\label{eq:DDDEb-pr33}
\end{equation}
which are two distinct real eigenvalues. 
We can conclude that $(M_-, M_-)$ is asymptotically stable and $(M_+, M_+)$ is a saddle. 
Indeed, this is because $M_- < \sqrt{k_3/k_1} < M_+$ holds by examining the quadratic function $F = F(W)$ derived from \eqref{eq:DDDEb-pr31}. 

The eigenvalues of $J(P_0)$ are $0$ and $-(1 + k_3) < 0$, and the corresponding eigenvectors (with $T$ as transpose) are
\[
{\mathbf{v}}_{1}=(1,1)^{T}, \quad 
{\mathbf{v}}_{2}=(1, -k_{3})^{T}.
\]
The dynamics near $P_0$ can be understood by translating the origin, applying normal form transformations, and applying the center manifold theorem. 
First, consider shifting the equilibrium $P_{0}$ to the origin:
\begin{equation}
\tilde{U}=U- \dfrac{k_2}{2k_1},\quad
\tilde{W}=W - \dfrac{k_2}{2k_1}.
\label{eq:DDDEb-pr34}
\end{equation}
The equations satisfied by $\tilde{U}$ and $\tilde{W}$ are
\begin{equation}
\begin{cases} 
\tilde{U}_{\tau} = -\tilde{U} + \tilde{W}, \\ 
\tilde{W}_{\tau} = k_3 \tilde{U} - k_3 \tilde{W} + k_2 \tilde{U}\tilde{W} - \dfrac{k_2}{2}\tilde{W}^2 + k_1 \tilde{U}\tilde{W}^2. 
\end{cases}
\label{eq:DDDEb-pr35}
\end{equation}
We examine the dynamics near the origin $(\tilde{U}, \tilde{W}) = (0, 0)$ for this system. 
Setting $P=(\mathbf{v}_{1},\mathbf{v}_{2})$, we obtain 
\begin{align*}
\left(\begin{array}{cc}
\tilde{U}_{\tau} \\
\tilde{W}_{\tau}
\end{array}
\right) &= \left(\begin{array}{cc}
-1 & 1 \\ k_{3} & -k_{3}
\end{array}
\right)\left(\begin{array}{cc}
\tilde{U} \\
\tilde{W}
\end{array}
\right)+\left(\begin{array}{cc}
0 \\
k_2 \tilde{U}\tilde{W} - \frac{k_2}{2}\tilde{W}^2 + k_1 \tilde{U}\tilde{W}^2
\end{array}
\right) \\
&= P\left(\begin{array}{cc}
0 & 0 \\
0 & -(1+k_{3})
\end{array}
\right)P^{-1}\left(\begin{array}{cc}
\tilde{U} \\
\tilde{W}
\end{array}
\right)+\left(\begin{array}{cc}
0 \\
k_2 \tilde{U}\tilde{W} - \frac{k_2}{2}\tilde{W}^2 + k_1 \tilde{U}\tilde{W}^2
\end{array}
\right).
\end{align*}
Let 
$\left(\begin{array}{cc}
\hat{U} \\
\hat{W}
\end{array}
\right)=P^{-1}\left(\begin{array}{cc}
\tilde{U} \\
\tilde{W}
\end{array}
\right)$.
Hence, we have 
\begin{equation}
\begin{cases} 
\hat{U}_{\tau}
= \dfrac{1}{1+k_3} \left[ \dfrac{k_2}{2}\hat{U}^2 + k_2 \hat{U}\hat{W} - \dfrac{k_2 k_3 (2+k_3)}{2}\hat{W}^2 \right] + \mathcal{O}(\vert{}(\hat{U}, \hat{W})\vert{}^3), 
\\ 
\hat{W}_{\tau}
= -(1+k_3)\hat{W} - \dfrac{1}{1+k_3} \left[ \dfrac{k_2}{2}\hat{U}^2 + k_2 \hat{U}\hat{W} - \dfrac{k_2 k_3 (2+k_3)}{2}\hat{W}^2 \right] + \mathcal{O}(\vert{}(\hat{U}, \hat{W})\vert{}^3). 
\end{cases}
\label{eq:DDDEb-pr36}
\end{equation}
Applying the center manifold theorem to \eqref{eq:DDDEb-pr36}, there exists a function $h(\hat{U})$ satisfying
\[
h(0)=\dfrac{dh}{d\hat{U}}(0)=0
\]
such that the center manifold of the origin is locally represented as $\{ (\hat{U}, \hat{W}) \mid \hat{W} = h(\hat{U}) \}$. 
The approximation of the center manifold is
\begin{equation}
\left\{ (\hat{U}(\tau), \hat{W}(\tau)) \mid \hat{W} = -\dfrac{k_2}{2(1+k_3)^2}\hat{U}^2 + \mathcal{O}(\hat{U}^3) \right\}.
\label{eq:DDDEb-pr37}
\end{equation}
Thus, the dynamics near the origin in \eqref{eq:DDDEb-pr36} is topologically equivalent to the dynamics of the following equation:
\begin{equation}
\hat{U}_{\tau} = \dfrac{k_2}{2(1+k_3)}\hat{U}^2 + \mathcal{O}(\hat{U}^3).
\label{eq:DDDEb-pr38}
\end{equation}
Consequently, the approximation of the center manifold near the origin for \eqref{eq:DDDEb-pr35} is
\begin{equation}
\left\{ (\tilde{U}(\tau), \tilde{W}(\tau)) \mid \tilde{W} = \tilde{U}+\dfrac{k_2}{2(1+k_3)}\tilde{U}^2 + \mathcal{O}(\tilde{U}^3) \right\}
\label{eq:DDDEb-pr39}
\end{equation}
and the dynamics near the origin for \eqref{eq:DDDEb-pr35} is topologically equivalent to:
\begin{equation}
\tilde{U}_{\tau} = \dfrac{k_2}{2(1+k_3)}\tilde{U}^2 + \mathcal{O}(\tilde{U}^3).
\label{eq:DDDEb-pr40}
\end{equation}
From the above discussion, the approximation of the center manifold near the equilibrium $P_0$ is given by
\begin{equation}
\left\{ (U(\tau), W(\tau)) \mid W = U + \frac{k_{2}}{2(1+k_{3})} \left( U - \frac{k_{2}}{2k_{1}} \right)^2 + \mathcal{O}\left( \left( U - \frac{k_{2}}{2k_{1}} \right)^3 \right)
 \right\},
\label{eq:DDDEb-pr41}
\end{equation}
and the dynamics around $P_0$ is topologically equivalent to:
\begin{equation}
U_{\tau}=\dfrac{k_2}{2(1+k_3)} \left( U - \frac{k_{2}}{2k_{1}} \right)^2 + \mathcal{O}\left( \left( U - \frac{k_{2}}{2k_{1}} \right)^3 \right).
\label{eq:DDDEb-pr42}
\end{equation}
These results show that a saddle-node bifurcation occurs when the equilibrium $P_0$ appears.

\subsection{Dynamics on the chart $\overline{U}_{2}$}
\label{sub:DDDEb-pr12}
Transformation
\[
U=\dfrac{\lambda_2}{\lambda_1}, \quad 
W=\dfrac{1}{\lambda_1}
\]
and time rescaling $d\eta/d\tau=\lambda_1^{-2}$ are induced to the following system in the $(\lambda_{1}, \lambda_{2})$-plane are
\begin{equation}
\begin{cases}
d\lambda_1/d\eta = -k_{1}\lambda_{1}\lambda_2 +k_2 \lambda_{1}^{2}- k_{3}\lambda_1^{3}, \\
d\lambda_2/d\eta = -\lambda_1^{2}\lambda_2 +\lambda_{1}^{2} -k_{1}\lambda_{2}^{2}+k_{2}\lambda_{1}\lambda_{2} -k_{3}\lambda_{1}^{2}\lambda_{2}.
\end{cases}
\label{eq:DDDEb-pr43}
\end{equation}
The equilibrium on $\{\lambda_{1}=0\}$ is $P_{1}: (\lambda_{1}, \lambda_{2})=(0,0)$, and the corresponding linearized matrix $J(P_{1})$ is 
\[
J(P_{1})=\left(\begin{array}{cc}
 0 & 0 \\ 0 & 0
\end{array}\right).
\]
Since $J(P_1)$ has double zero eigenvalues, the dynamics near $P_1$ can be understood using the blow-up transformation
\[
\lambda_{1}=r\bar{\lambda}_{1}, \quad \lambda_{2}=r\bar{\lambda}_{2}.
\]

\subsubsection{Dynamics on the chart $\{\bar{\lambda}_{1}=1\}$}
\label{sub:DDDEb-pr121}
By using the transformation $\lambda_1 = r$, $\lambda_2 = r \bar{\lambda}_2$ and the time scale transformation $ds/d\eta = r$, we obtain
\begin{equation}
\begin{cases}
dr/ds= -k_{1}r\bar{\lambda}_{2}+k_{2}r -k_{3}r^{2}, \\
d\bar{\lambda}_{2}/ds= -r\bar{\lambda}_{2}+1.
\end{cases}
\label{eq:DDDEb-pr44}
\end{equation}
No equilibrium exists on $\{r = 0\}$, and $d\bar{\lambda}_2/ds > 0$ holds on $\{r = 0\}$.

\subsubsection{Dynamics on the chart $\{\bar{\lambda}_{2}=1\}$}
\label{sub:DDDEb-pr122}
By introducing $\lambda_1 = r \bar{\lambda}_1$, $\lambda_2 = r$ and $ds/d\eta = r$, we obtain
\begin{equation}
\begin{cases}
dr/ds=-r^{2}\bar{\lambda}_{1}^{2}+r\bar{\lambda}_{1}^{2} -k_{1}r+k_{2}r\bar{\lambda}_{1} -k_{3}r^{2}\bar{\lambda}_{1}^{2}, \\
d\bar{\lambda}_{1}/ds = r\bar{\lambda}_{1}^{3} - \bar{\lambda}_{1}^{3}.
\end{cases}
\label{eq:DDDEb-pr45}
\end{equation}
The equilibrium on $\{r = 0\}$ is $P_2: (r, \bar{\lambda}_1) = (0, 0)$, and the linearized matrix $J(P_2)$ at this equilibrium is
\[
J(P_{2})=\left(\begin{array}{cc}
-k_{1} & 0 \\ 0 & 0 
\end{array}\right).
\]
The dynamics near $P_2$ follows the same argument as in Subsection  \ref{sub:DDDEb-pre72}; the invariant set $\{r = 0\}$ is one representation of the center manifold, and the dynamics on this center manifold is given by \eqref{eq:DDDEb-pr21}.

\subsection{Dynamics on the chart $\overline{U}_{1}$}
\label{sub:DDDEb-pr13}
We introduce the ransformation
\[
U=\dfrac{1}{\lambda_{1}}, \quad W=\dfrac{\lambda_2}{\lambda_1}
\]
and time scale transformation $d\eta/d\tau =\lambda_{1}^{-2}$.
We then obtain
\begin{equation}
\begin{cases}
d\lambda_1/d\eta
= \lambda_{1}^{3}-\lambda_{1}^{3}\lambda_{2}, \\ 
d\lambda_2/d\eta
= k_{1}\lambda_{2}^{2}-k_{2}\lambda_{1}\lambda_{2}^{2}+(k_{3}+1)\lambda_{1}^{2}\lambda_{2}-\lambda_{1}^{2}\lambda_{2}^{2}.
\end{cases}
\label{eq:DDDEb-pr46}
\end{equation}
The equilibrium on $\{\lambda_1 = 0\}$ is $P_3: (\lambda_1, \lambda_2) = (0, 0)$. Since its linearized matrix easily turns out to have double zero eigenvalues, we can understand the dynamics near $P_3$ via the blow-up transformation:
\[
\lambda_{1}=r\bar{\lambda}_{1}, \quad 
\lambda_{2}=r^{2}\bar{\lambda}_{2}.
\]

\subsubsection{Dynamics on the chart $\{\bar{\lambda}_{1}=1\}$}
\label{sub:DDDEb-pr131}
Using $\lambda_{1}=r$, $\lambda_{2}=r^{2}\bar{\lambda}_{2}$, and the time scale transformation $ds/d\eta = r^2$, we obtain
\begin{equation}
\begin{cases}
dr/ds= r-r^{3}\bar{\lambda}_{2}, \\
d\bar{\lambda}_{2}/ds=(k_{3}-1)\bar{\lambda}_{2}+ r^{2}\bar{\lambda}_{2}^{2}+ k_{1} \bar{\lambda}_{2}^{2}- k_{2} r\bar{\lambda}_{2}^{2}.
\end{cases}
\label{eq:DDDEb-pr47}
\end{equation}
As equilibria on $\{r = 0, \bar{\lambda}_2 \ge 0\}$, $P_{4}: (r, \bar{\lambda}_{2})=(0,0)$ exists, and when $0 < k_{3} < 1$, $P_{5}: (r, \bar{\lambda}_2) = (0, [1 - k_{3}]/k_{1})$ exists. 
The linearized matrices at these equilibria are
\[
J(P_{4})=\left(\begin{array}{cc}
1 & 0 \\ 0 & k_{3}-1
\end{array}\right), \quad
J(P_{5})=\left(\begin{array}{cc}
1 & 0 \\ -k_{2}[1-k_3]^{2}/k_{1}^{2} & 1-k_{3}
\end{array}\right).
\]
Therefore, we conclude the following results:
\begin{enumerate}
\item[(i)] When $0<k_{3}<1$, $P_{4}$ is a saddle and $P_{5}$ is unstable.
\item[(ii)] When $k_{3}=1$, $P_{4}$ is not hyperbolic since the corresponding linearized matrix has a zero eigenvalue.
However, $dr/ds = 0$ and $d\bar{\lambda}_2/ds > 0$ hold near $r = 0$.
We can conclude that no trajectory converges to $P_4$.  
\item[(iii)] When $k_3 > 1$, $P_4$ is an unstable equilibrium.
\end{enumerate}

\subsubsection{Dynamics on the chart $\{\bar{\lambda}_{2}=1\}$}
\label{sub:DDDEb-pr132}
By using the transformation $\lambda_1 = r \bar{\lambda}_1$, $\lambda_2 = r^2$ and the time scale transformation $ds/d\eta = r^2$, we obtain
\begin{equation}
\begin{cases}
dr/ds= 2^{-1}k_{1} r - 2^{-1}k_{2} r^{2}\bar{\lambda}_{1} + 2^{-1}(k_{3}+1) r\bar{\lambda}_{1}^{2}-2^{-1}r^{3}\bar{\lambda}_{1}^{2}, \\
d\bar{\lambda}_{1}/ds = -2^{-1}k_{1} \bar{\lambda}_{1}+2^{-1}k_{2}r \bar{\lambda}_{1}^{2} +2^{-1}(1-k_{3}) \bar{\lambda}_{1}^{3}-2^{-1}r^{2}\bar{\lambda}_{1}^{3}.
\end{cases}
\label{eq:DDDEb-pr48}
\end{equation}
The equilibrium on $\{r = 0, \bar{\lambda}_1 \ge 0\}$ is $P_{6}: (r, \bar{\lambda}_1) = (0, 0)$, and when $0 < k_{3} < 1$, $P_{7}: (r, \bar{\lambda}_1) = (0, \sqrt{k_{1} / [1 - k_{3}]})$ exists. 
The linearized matrices near these points are
\[
J(P_{6})=\left(\begin{array}{cc}
2^{-1}k_{1} & 0 \\ 0 & -2^{-1}k_{1}
\end{array}\right), \quad
J(P_{7})= \left(\begin{array}{cc}
k_{1}/[1-k_{3}] & 0 \\ 2^{-1}k_{2}k_{1}/[1-k_{3}] & k_1
\end{array}\right).
\]
Therefore, we conclude that followings hold:
\begin{enumerate}
\item[(i)] When $0<k_{3}<1$, $P_{6}$ is a saddle and $P_{7}$ is an unstable equilibrium.
\item[(ii)] When $k_{3}=1$, $P_{6}$ is a saddle.
\item[(iii)] When $k_{3}>1$, $P_{6}$ is a saddle.
\end{enumerate}
Combining the above discussions, we see that no trajectory is attracted to the origin on the local chart $\overline{U}_1$.

\subsection{Proof of Theorem \ref{th:DDDEb-mr3}}
\label{sub:DDDEb-pr14}
\begin{figure}[h]
\begin{center}
\includegraphics[scale=0.21]{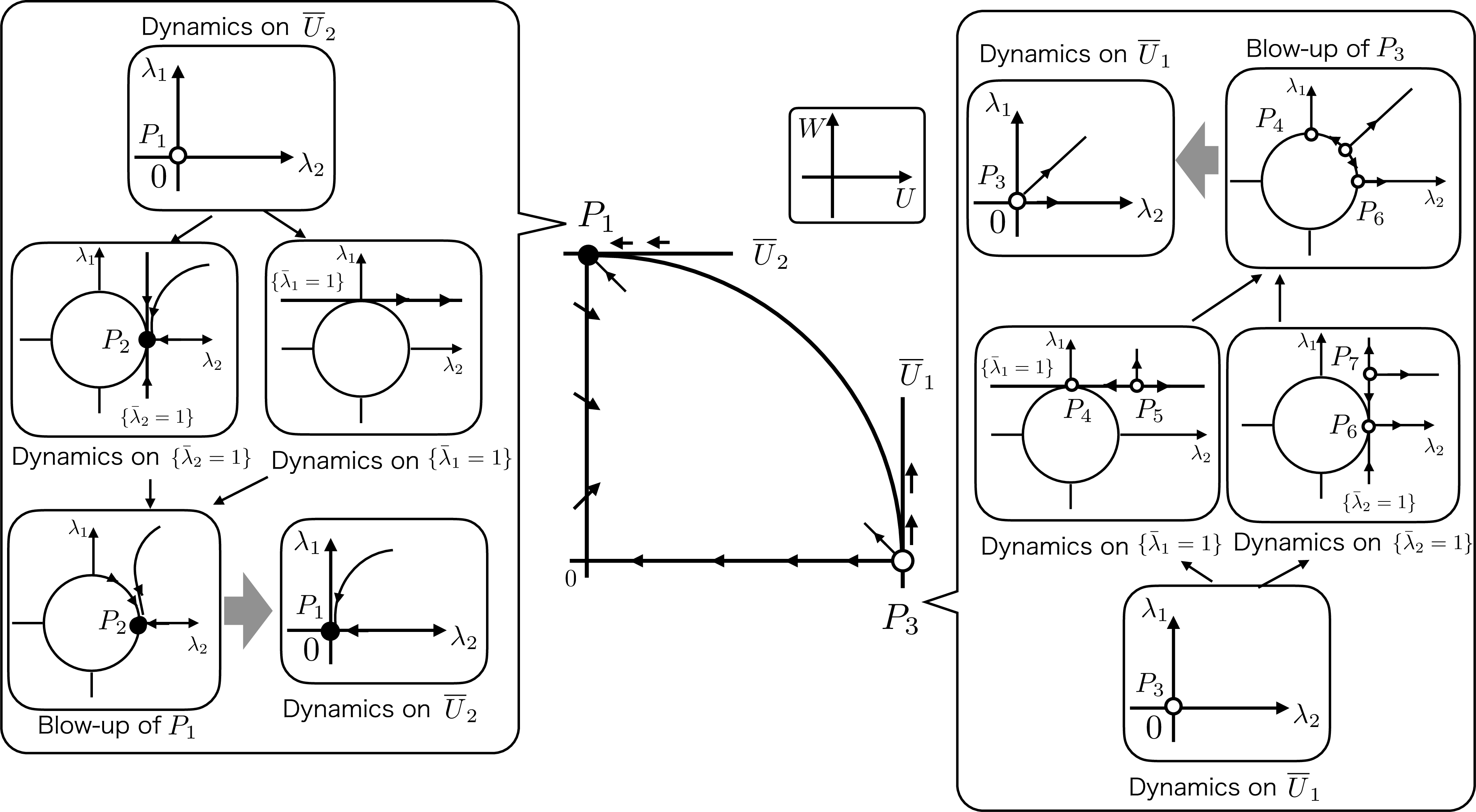}
\caption{Schematic pictures of the dynamics on a Poincar\'e-type disk in \eqref{eq:DDDEb-pr28}.}
\label{fig:DDDEb5}
\end{center}
\end{figure}

By combining the dynamics on local charts $\overline{U}_1$ and $\overline{U}_2$ for \eqref{eq:DDDEb-pr28}, we obtain the dynamics including infinity as shown in Figure \ref{fig:DDDEb5}. 

\begin{figure}[h]
\begin{center}
\includegraphics[scale=0.25]{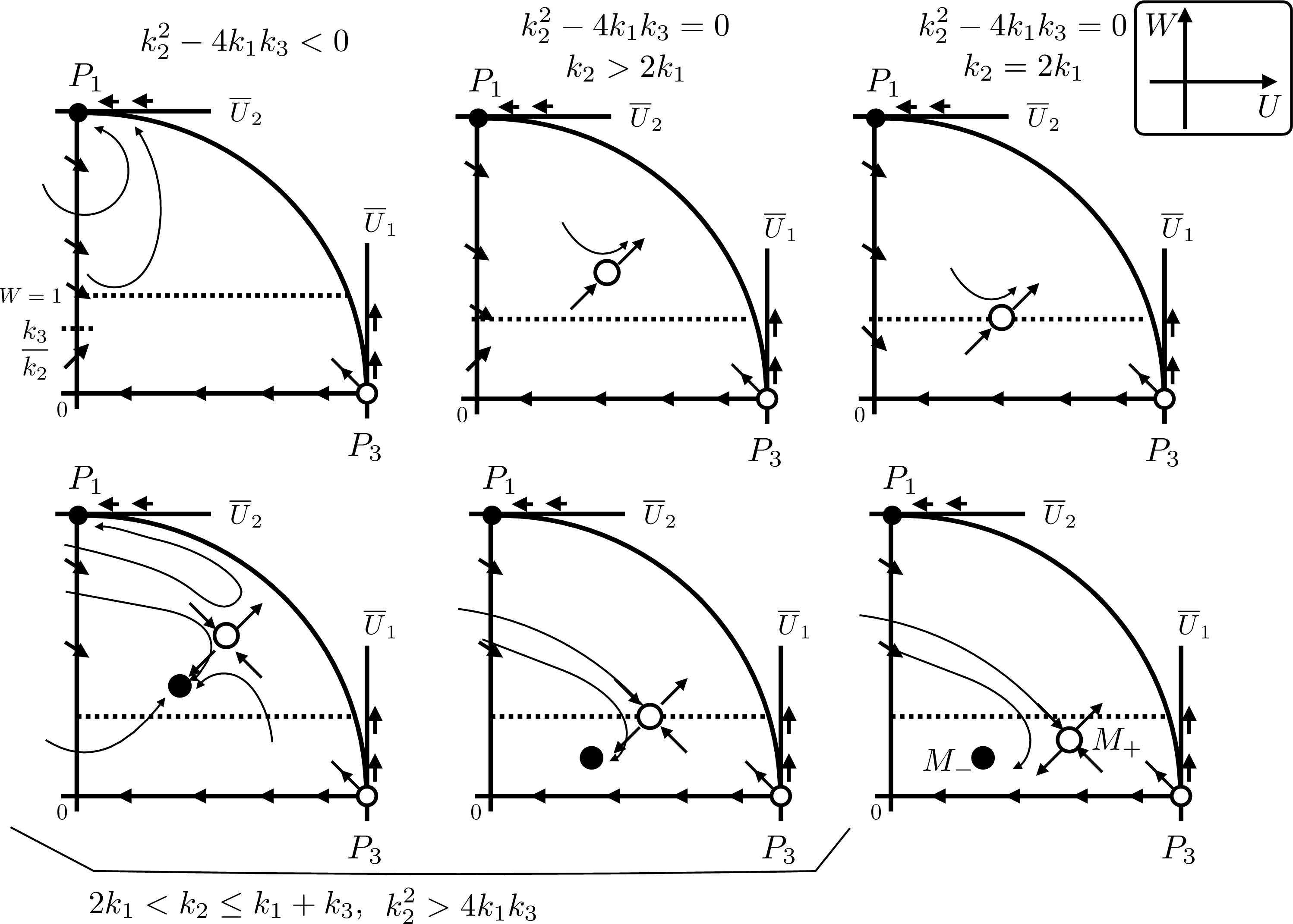}
\caption{Schematic pictures of the dynamics near the finite eqauilibrium in \eqref{eq:DDDEb-pr28}.}
\label{fig:DDDEb6}
\end{center}
\end{figure}

Combining these with the dynamics near finite equilibria yields schematic pictures as shown in Figure \ref{fig:DDDEb6}. 
When $k_{2}^{2}-4k_{1}k_{3}<0$, no finite equilibrium exists in $\Phi_4$. Moreover, since $k_{2}>k_{3}>0$ holds, $k_{3}/k_{2}<1$. 
This implies that $W_\tau < 0$ holds in the region on the $W$-axis satisfying $W > 1$ in \eqref{eq:DDDEb-pr28}. 
Thus, arguments similar to Sections \ref{sec:DDDEb-pr1} and \ref{sec:DDDEb-pr2} show that trajectories starting from sufficiently large initial values $(U(0), W(0))$ head toward the equilibrium $P_1$. 
However, since $W_\tau < 0$ on $\{(U, W) \mid U<k_{1}^{-1}[k_{2}-k_{3}],\,\, W=1\}$, if an initial value $(U(0), W(0))$ is chosen close to $(U, W) = (0, 1)$ on the $UW$-plane, the trajectory starting from there crosses the line $W = 1$ and exits the region $\Phi_4$. 
That is, depending on the choice of initial values $(U(0), W(0))$, there exists a threshold curve on the phase plane governing whether trajectories exit $\Phi_4$ or remain in $\Phi_4$ to eventually head toward $P_1$. 

When $k_{2}^{2}-4k_{1}k_{3}=0$, even if the equilibrium $P_0$ lies on $W = 1$, initial values in $\Phi_4$ attracted to $P_1$ can be chosen. 
On the other hand, when $k_{2}<2k_{1}$, the equilibrium $P_0$ lies below $W = 1$. 
In this case, there exists a trajectory starting from a positive point $(U, W) = (0, W_c)$ on the $W$-axis and converging to $P_0$. 
Initial values $(U(0), W(0))$ can be taken below this trajectory, implying a threshold phenomenon regarding the occurrence/non-occurrence of blow-up depending on initial conditions and parameters.

Now consider $k_{2}^{2}-4k_{1}k_{3}>0$. 
Regardless of the value of $M_+$, there exists a trajectory starting from a positive point on the $W$-axis and converging to $P_+$. 
On the $UW$-plane, trajectories starting above the curve corresponding to this trajectory converge to $P_1$, representing blow-up solutions. 
On the other hand, trajectories starting below this curve converge to $P_-$, which means they do not correspond to blow-up solutions.

It turns out that due to the self-inhibitory term, finite-time singularities do not occur if and only if the initial state $u(0)$ and accumulated memory $v(0)$ at $t = 0$ are sufficiently small. 
In particular, when $k_{2}^{2}-4k_{1}k_{3}>0$ and $M_- > 1$, even if $u(0)$ is large, the solution exists globally without blowing up in finite time as long as $w(0)$ is small. 
This implies that if the accumulated memory $v(0)$ at $t = 0$ is small, the effect of the self-inhibitory term works effectively even if $u(0)$ is large.

What remains to be shown is that finite-time blow-up occurs and to derive the blow-up rates \eqref{eq:DDDEb-mr3}. 
From \eqref{eq:DDDEb-pr45}, the leading part of the dynamics near $P_2$ in Subsection \ref{sub:DDDEb-pr122} is $dr/ds \sim -k_1 r$, and the argument proceeds by replacing $\alpha$ in Subsection \ref{sub:DDDEb-pr9} with $k_1$. 
Consequently, we obtain \eqref{eq:DDDEb-mr3} as the blow-up rate.

\section{Discussion}
\label{sec:DDDEb-d}
In this paper, we clarify the existence, blow-up rates, and threshold phenomena of finite-time blow-up solutions for distributed delay differential equations with a gamma distribution kernel by specifying concrete nonlinearities. 
In particular, the nonlinearities considered in this paper serve as typical examples of finite-time singularities in solutions of ODEs, representing essential nonlinear terms that should be examined as a first step when discussing finite-time singularities.

In equation \eqref{eq:DDDEb-int1}, we consider
\begin{equation}
\dot{u}(t) = f\left( u, \int_{-t_{c}}^{t} \alpha e^{-\beta(t-s)} g(u(s)) \, ds \right)
\label{eq:DDDEb-f1}
\end{equation}
with the positive finite value $t_{c}$.
Let $v_{t_c}(t)$ be defined as:
\begin{equation}
v_{t_{c}}(t)=\int_{-t_{c}}^{t} \alpha e^{-\beta(t-s)} g(u(s)) \, ds.
\label{eq:DDDEb-f2}
\end{equation}
By assuming that the initial function is positive and continuous, $v_{t_c}(0)$ at $t = 0$ takes a positive finite value. 
In this case, the Linear Chain Trick remains effective, yielding system \eqref{eq:DDDEb-pr2}. 
Since the conclusions in this paper depend strictly on the structure of the equations, the existence of blow-up solutions and their blow-up rates are identical to Theorems \ref{th:DDDEb-mr1}, \ref{th:DDDEb-mr2}, and \ref{th:DDDEb-mr3} in Section \ref{sec:DDDEb-mr}. 
However, since $v_{t_c}(0) < v(0)$ holds, care must be taken when precisely analyzing the relationship between the initial function and blow-up solutions

In this section, for Theorems \ref{th:DDDEb-mr1}, \ref{th:DDDEb-mr2}, and \ref{th:DDDEb-mr3} which adopt typical non-delayed ODE examples of finite-time singularities \eqref{eq:DDDEb-int3} and \eqref{eq:DDDEb-int6}, we derive the corresponding non-delayed and memory-free differential equations by taking the limit $\beta \rightarrow +\infty$. 
By comparing the behaviors with non-delayed differential equations lacking memory effects, we can discuss the influence of memory accumulation and forgetting, as well as finite-time singularities induced by time delays. 
In this argument, two distinct developments arise depending on whether $\alpha = \beta$ or $\alpha \neq \beta$, both of which are detailed below.

From the definition of $v(t)$ in \eqref{eq:DDDEb-int2}, setting $x = t - s$ and $y = \beta(t - s)$ yields:
\begin{align*}
v(t) 
&= \frac{\alpha}{\beta} \int_{-\infty}^{t} \beta e^{-\beta(t-s)} g(u(s)) \, ds
= \frac{\alpha}{\beta} \int_{0}^{\infty} \beta e^{-\beta x} g(u(t-x)) \, dx
\\
&= \frac{\alpha}{\beta} \int_{0}^{\infty} e^{-y} g\left(u\left(t - \frac{y}{\beta}\right)\right) \, dy. 
\end{align*}
Fix an arbitrary time $t \in (0, T)$. 
Since the initial function is bounded from \eqref{eq:DDDEb-ini1}, there exists a constant $K_t > 0$ depending on $t$ such that for all $y \ge 0$ and $\beta > 0$:
\begin{equation}
\left| g\left(u\left(t - \frac{y}{\beta}\right)\right) \right| \le K_{t}.
\label{eq:DDDeb-d01}
\end{equation}
Since $g(u)$ specified in (F1), (F2), and (F3) is smooth, it holds that:
\begin{equation}
\left| e^{-y} g\left(u\left(t - \frac{y}{\beta}\right)\right) \right| \le K_{t}e^{-y}. 
\label{eq:DDDeb-d02}
\end{equation}
Therefore, for each $t \in (0, T)$, applying Lebesgue's dominated convergence theorem gives as follows:
\begin{align*}
\lim_{\beta\to +\infty}\int_{0}^{\infty} e^{-y} g\left(u\left(t - \frac{y}{\beta}\right)\right) \, dy 
&=
\int_{0}^{\infty} \lim_{\beta\to +\infty} e^{-y} g\left(u\left(t - \frac{y}{\beta}\right)\right) \, dy \\
&= g(u(t)) \int_{0}^{\infty} e^{-y} \, dy \\
&= g(u(t)). 
\end{align*}
Then we have 
\begin{equation}
\lim_{\beta\to +\infty} v(t)
=\begin{cases}
g(u(t)) \quad & (\alpha=\beta), \\
0 \quad & (\alpha\neq \beta).
\end{cases}
\label{eq:DDDEb-d1}
\end{equation}
Note that in the case $\alpha \neq \beta$, $\alpha$ is held fixed while $\beta \rightarrow +\infty$, meaning $\alpha$ is independent of $\beta$.  

By setting $\alpha = \beta$ and taking $\beta \rightarrow +\infty$ in \eqref{eq:DDDEb-d1}, the non-delayed differential equation without memory effects corresponding to the delay differential equation \eqref{eq:DDDEb-int1} is derived as:
\begin{equation}
\dot{u}=f(u(t), g(u(t))).
\label{eq:DDDEb-d2}
\end{equation}
Specifically, when (F1) is imposed on \eqref{eq:DDDEb-int1}, the corresponding non-delayed (memery-free) differential equation becomes:
\begin{equation}
\dot{u}(t)= \{ g(u(t)) \}^{p}= u(t)^{p}.
\label{eq:DDDEb-d3}
\end{equation}
When (F2) is imposed, it becomes:
\begin{equation}
\dot{u}(t)= \dfrac{1}{1-g(u(t))}=\dfrac{1}{1-u(t)}.
\label{eq:DDDEb-d4}
\end{equation}
Regarding the solution behavior of the corresponding non-delayed differential equations, (F1) and (F2) correspond to \eqref{eq:DDDEb-int3} and \eqref{eq:DDDEb-int6}, respectively.
Under the case for $\alpha = \beta$, the memory structure is introduced as accumulation via $v(t)$.
In the case of (F1), the blow-up rate in \eqref{eq:DDDEb-int5} is accelerated to \eqref{eq:DDDEb-mr1} due to the effect of $\beta$. 
This suggests that when $\beta > 0$ is finite, the memory accumulation $v(t)$ acts as the dominant factor of the blow-up structure, providing positive feedback as a blow-up booster.
On the other hand, in the case of (F2), a dramatic qualitative shift is observed: while the non-delayed equation undergoes quenching, the memory structure and $v(t)$ induce blow-up, with the logarithmic rate given by \eqref{eq:DDDEb-mr3}.
For (F2), this represents a case where blow-up is induced by the effect of memory accumulation and forgetting, which in this paper is termed ``memory-induced blow-up''.
In addition, since the dynamics undergo a dramatic qualitative transition due to memory effects, this phenomenon can also be described as a ``qualitative transition induced by memory''.
Similarly for (F3), the non-delayed differential equation is as follows:
\begin{equation}
\dot{u}(t)=\dfrac{1}{1-[u(t)-a]}-\gamma u(t),
\label{eq:DDDEb-d5}
\end{equation}
which implies that memory effects acting as time delays induce memory-induced blow-up, and that memory induces a threshold governing global existence in time. See also Table \ref{tb:DDDEb1}.

We also consider the limiting non-delayed ODE derived by taking $\beta \rightarrow +\infty$ while holding $\alpha$ fixed, with $\alpha$ and $\beta$ treated as independent parameters. 
From \eqref{eq:DDDEb-d1}, the limiting ordinary differential equation corresponding to \eqref{eq:DDDEb-int1} is given by:
\begin{equation}
\dot{u}(t)=\hat{f}(u):=f(u, 0).
\label{eq:DDDEb-d6}
\end{equation}
Under assumption (F1), the resulting limiting equation is $\dot{u} = 0$, whose solution is a constant function.
Thus, memory effects turn a static state into blow-up, demonstrating the occurrence of memory-induced blow-up. 
Similarly, for (F2), the non-delayed equation is $\dot{u} = 1$, and for (F3), it is $\dot{u} = 1 - \gamma u(t)$. 
See also Table \ref{tb:DDDEb1}. 
Consequently, in (F2), memory effects serve to transition a globally existing solution (infinite-time divergence) into a finite-time blow-up solution.
In (F3), memory destroys system stability and induces blow-up phenomena.  

\begin{table}[htbp]
\centering
\begin{tabular}{c||c|c}
& $\alpha=\beta$ & $\alpha\neq \beta$ \\ \hline\hline
(F1) & $\dot{u}(t)=u(t)^{p}$ (Rate Acceleration) & $\dot{u}(t)=0$ (Static $\to$ Blow-up) 
\\[11pt] \hline
(F2) & $\dot{u}(t)=\dfrac{1}{1-u(t)}$ (Quench $\to$ Blow-up) & $\dot{u}(t)=1$ (Global $\to$ Blow-up) 
\\[11pt] \hline
(F3) & $\dot{u}(t)=\dfrac{1}{1-[u(t)-a]}-\gamma u(t)$ (Threshold) & $\dot{u}(t)=1-\gamma u(t)$ (Stable $\to$ Threshold) 
\\ \hline
\end{tabular}
\caption{Qualitative comparison between the dynamics of non-delayed limiting ODEs (derived formally as $\beta \to +\infty$) and the corresponding distributed delay differential equations}
\label{tb:DDDEb1}
\end{table}

\section*{Acknowledgments}
This work was partly supported by the JSPS KAKENHI (Grant No. 25K17306).
The author used Gemini (google) for language editing and assistance with manuscript preparation. 
All AI-generated content was reviewed, edited, and approved by the author, who takes full responsibility for the integrity of the content.


\end{document}